\documentclass[12pt]{amsart}
\usepackage{amsbsy}
\usepackage{graphicx,epsfig,subfigure,psfrag}
\begin{document}

\newtheorem{theorem}{Theorem}
\newtheorem{proposition}{Proposition}
\newtheorem{lemma}{Lemma}
\newtheorem{corollary}{Corollary}
\newtheorem{definition}{Definition}
\newtheorem{remark}{Remark}
\newcommand{\tex}{\textstyle}
\numberwithin{equation}{section} \numberwithin{theorem}{section}
\numberwithin{proposition}{section} \numberwithin{lemma}{section}
\numberwithin{corollary}{section}
\numberwithin{definition}{section} \numberwithin{remark}{section}
\newcommand{\ren}{\mathbb{R}^N}
\newcommand{\re}{\mathbb{R}}
\newcommand{\n}{\nabla}
\newcommand{\iy}{\infty}
\newcommand{\pa}{\partial}
\newcommand{\fp}{\noindent}
\newcommand{\ms}{\medskip\vskip-.1cm}
\newcommand{\mpb}{\medskip}
\newcommand{\BB}{{\bf B}}
\newcommand{\Am}{{\bf A}_{2m}}
\renewcommand{\a}{\alpha}
\renewcommand{\b}{\beta}
\newcommand{\g}{\gamma}
\newcommand{\G}{\Gamma}
\renewcommand{\d}{\delta}
\newcommand{\D}{\Delta}
\newcommand{\e}{\varepsilon}
\newcommand{\var}{\varphi}
\renewcommand{\l}{\lambda}
\renewcommand{\o}{\omega}
\renewcommand{\O}{\Omega}
\newcommand{\s}{\sigma}
\renewcommand{\t}{\tau}
\renewcommand{\th}{\theta}
\newcommand{\z}{\zeta}
\newcommand{\wx}{\widetilde x}
\newcommand{\wt}{\widetilde t}
\newcommand{\noi}{\noindent}
\newcommand{\uu}{{\bf u}}
\newcommand{\UU}{{\bf U}}
\newcommand{\VV}{{\bf V}}
\newcommand{\ww}{{\bf w}}
\newcommand{\vv}{{\bf v}}
\newcommand{\WW}{{\bf W}}
\newcommand{\hh}{{\bf h}}
\newcommand{\di}{{\rm div}\,}
\newcommand{\inA}{\quad \mbox{in} \quad \ren \times \re_+}
\newcommand{\inB}{\quad \mbox{in} \quad}
\newcommand{\inC}{\quad \mbox{in} \quad \re \times \re_+}
\newcommand{\inD}{\quad \mbox{in} \quad \re}
\newcommand{\forA}{\quad \mbox{for} \quad}
\newcommand{\whereA}{,\quad \mbox{where} \quad}
\newcommand{\asA}{\quad \mbox{as} \quad}
\newcommand{\andA}{\quad \mbox{and} \quad}
\newcommand{\withA}{,\quad \mbox{with} \quad}
\newcommand{\orA}{,\quad \mbox{or} \quad}
\newcommand{\ssk}{\smallskip}
\newcommand{\LongA}{\quad \Longrightarrow \quad}
\def\com#1{\fbox{\parbox{6in}{\texttt{#1}}}}
\def\N{{\mathbb N}}
\def\A{{\cal A}}
\newcommand{\de}{\,d}
\newcommand{\eps}{\varepsilon}
\newcommand{\be}{\begin{equation}}
\newcommand{\ee}{\end{equation}}
\newcommand{\spt}{{\mbox spt}}
\newcommand{\ind}{{\mbox ind}}
\newcommand{\supp}{{\mbox supp}}
\newcommand{\dip}{\displaystyle}
\newcommand{\prt}{\partial}
\renewcommand{\theequation}{\thesection.\arabic{equation}}
\renewcommand{\baselinestretch}{1.2}

\title
{\bf On global solutions and blow-up for\\
Kuramoto--Sivashinsky-type models,\\
 and   well-posed Burnett
  equations}


\author{
V.A.~Galaktionov, E. Mitidieri,  and S.I.~Pohozaev}

\address{Department of Mathematical Sciences, University of Bath,
 Bath BA2 7AY, UK}
\email{vag@maths.bath.ac.uk}

\address{Dipartimento di Matematica e Informatica,
Universit\`a  di Trieste, Via  Valerio 12, 34127 Trieste, ITALY}
 \email{mitidier@units.it}

\address{Steklov Mathematical Institute,
 Gubkina St. 8, 119991 Moscow, RUSSIA}
\email{pokhozhaev@mi.ras.ru}



  \keywords{Higher-order semilinear parabolic equations, global existence, uniform bounds,
  nonexistence, blow-up, Navier-Stokes and  Burnett equations.
  {\bf To appear in:} Nonl. Anal.}
 \date{\today}




 \begin{abstract}

 The  initial boundary-value problem (IBVP) and the Cauchy problem
for the Kuramoto--Sivashinsky equation
 $$
 \tex{
v_t+v_{xxxx}+v_{xx}=\frac{1}{2}\, (v^2)_x
 }
 $$
 and other related $2m$th-order semilinear parabolic partial differential equations
 in one dimension and in $\ren$ are considered.
Global existence and blow-up as well as $L^\infty$-bounds are
reviewed by using:

 (i) classic tools of interpolation theory and Galerkin methods,

(ii) eigenfunction and nonlinear capacity methods,

 (iii) Henry's version of weighted Gronwall's inequalities,

(vi) two types of scaling (blow-up) arguments.

 \noi  For the IBVPs,
  existence of global solutions is proved for both Dirichlet and
``Navier" boundary conditions.
  For some related $2m$th-order PDEs in $\ren \times \re_+$,
uniform boundedness of global solutions of  the Cauchy problem are
established.

As another related application,  the well-posed Burnett-type
equations
 $$
 \vv_t +(\vv \cdot \n)\vv=- \n p - (-\D )^m \vv, \quad \di \vv=0 \inA, \quad m \ge
 1,
 $$
  are studied. For $m=1$, these are the classic Navier-Stokes
equations.
 As a simple illustration, it is shown  that a uniform $L^p(\ren)$-bound on locally sufficiently
  smooth $\vv(x,t)$
 for $p > \frac N{2m-1}$
  implies a uniform
  $L^\infty(\ren)$-bound, hence the solutions do not blow-up.
  For $m=1$, $N=3$, this gives $p > 3$ reflecting the famous
Leray--Prodi--Serrin--Ladyzhenskaya regularity results ($L^{p,q}$
criteria), and re-derives  Kato's class of unique mild solutions
in $\ren$. Bounded classic $L^2$-solutions are shown to  exist for
 $N < 2(2m-1)$.



\end{abstract}

\maketitle


\section{Introduction and motivation}
 \label{SInt}

 The role of the {\em Kuramoto--Sivashinsky
equation}
\begin{equation}
\label{e2.100}
 \tex{
v_t+v_{xxxx}+v_{xx}=\frac{1}{2}\, (v^2)_x
 }
\end{equation}
is well known in contemporary nonlinear mechanics and physics.

 This equation
arises as a model in hydrodynamics (a thin film flow down an
inclined plane in the presence of an electric field), in
combustion theory (propagation of flame fronts), phase turbulence
and plasma physics, as well as  one model for spatio-temporal chaos and in
many other physical phenomena; see \cite{9} for a nice short
review of applications with key original references.

First results on global existence of classical solutions of
(\ref{e2.100}) go back to the 1970s (cf. \cite{Lin74}) and the
1980s; we refer to \cite{1, Biag96, 2, 3, 5, 4, Gr03, Kai06, 6,
Kuk05, Lar04, Sell92, Tad86, 9,  10}, where further important references
can be found.
 A large part of  previous studie of  (\ref{e2.100}) was
 devoted to periodic problems. See \cite{2, Kuk05} for important references
 from the 1980s and also to some classes of particular solutions and the local behavior of
solutions, where very interesting results (for instance,  chaotic
behaviour, estimates of dimension and structure of attractors,
 bifurcation theorems, etc.) have been obtained. Other papers with deep
results are devoted to existence of  periodic solutions and traveling wave
solutions with complicated dynamics. Numerous contributions are dealing
with dynamical analysis of the one dimensional Kuramoto--Sivashinsky equation.





The present paper is devoted to a review of approaches that lead
to global existence and blow-up for the {\em Kuramoto--Sivashinsky
equation}  and modified versions of it  in one dimension and on $\ren$, the latter being much less developed in the literature.

One of the key mathematical features of the KS-type PDEs is that
an {\em a priori} $L^2$-bound of solutions $v(x,t)$ of the form
 \be
 \label{aa1}
  \tex{
  \|v(t)\|_2 = \displaystyle{\Big(\int |v(x,t)|^2 \, {\mathrm d}x}\Big)^{\frac12}
    \le \|v_0\|_2 \, {\mathrm e}^{\frac t4} \forA t \ge 0,
  }
  \ee
 is
straightforward, while the main question is how to use this information  to
get stronger estimates in Sobolev spaces and eventually in
$L^\infty$.

 To this end, we shall use some known methods and
compare their strength with the new techniques to be reviewed and developed throughout this paper:

\ssk

{\bf (I)}  In Section \ref{SGl} we employ classical interpolation and Galerkin methods for study global
existence
of the initial-boundary value problems (IBVP) in one dimension
with Dirichlet boundary conditions and with ``Navier" ones;

\ssk

{\bf (II)} Eigenfunction technique and nonlinear capacity
method \cite{MitPoh}  are used in Section \ref{SBlowup} to prove blow-up in
finite time of solutions of these IBVPs with non-standard boundary
conditions;

\ssk

{\bf (III)} In Section \ref{SGlH},  Henry's version of weighted Gronwall's inequalities are used
prove global existence for the Cauchy problem  in $\ren\times
\re_+$ for modified Kuramoto-Sivashinsky equations;

\ssk

{\bf (IV)} In Section \ref{SScal}, we use two types of scaling
blow-up techniques for global existence for the Cauchy problem and
IBVPs.

\ssk

 Our analysis embraces  a number of $2m$ th-order
$N$-dimensional Kuramoto-Sivashinsky  type models posed in $\ren
\times \re_+$.

\ssk

In general, as customary in PDE theory,
 the following approach  is classical:
 \be
  \label{St1}
   \fbox{$
  \mbox{local existence} + \mbox{a global {\em a priori} bound}
 \LongA \mbox{global existence}.
 $}
  \ee
 More precisely, in order to get (\ref{St1}), we use the following
 intermediate blow-up step:
\be
  \label{St2}
  \mbox{local existence} + \mbox{a global {\em a priori} bound}
 \LongA \mbox{no finite-time blow-up}.
  \ee
Here scaling arguments to prevent blow-up are crucial.

In this context, the problem of global solvability takes the
following obvious negation form:
\be
  \label{St3}
   \fbox{$
  \mbox{global solvability} = \not\exists \,\, \mbox{finite-time blow-up of local smooth solutions}.
 $}
  \ee
  Of course, here we understand the {\em a priori} estimates in the corresponding function
   spaces.
This proclamation is also  relevant to the Navier--Stokes
equations in $\re^3$, which have a long mysterious history of
uniquenes/non-uniqueness and blow-up singularity open problems.



As rather unexpected but  related application, in Section
\ref{SNS}, we consider the {\em Navier--Stokes equations}
  \be
  \label{NS1}
 \vv_t +(\vv \cdot \n)\vv=- \n p + \D \vv, \quad \di \vv=0 \inA,
  \ee
 with bounded integrable divergence-free data $\vv_0$. It is worth noting that the
 {\em convective term} in the first equation in (\ref{NS1}) for the
 velocity field has indeed a nonlinear dispersion mathematical
 nature as in the Kuramoto-Sivashinsky equation (\ref{e2.100}). Therefore, it is rather  natural
 to include the model (\ref{NS1}) in the present context of
 semilinear KS-type equations.

\ssk

\noi\underline{\em Blow-up self-similar singularities with finite
energy do not exist}. The idea that the classical fundamental
problem of unique solvability of (\ref{NS1}) in $\re^3$ is
associated with existence or nonexistence of certain {\em blow-up
singularities} as $t \to
T^-$, goes back to Th.~von~K\'arm\'an; see 
\cite{vonK21}. Later on, in 1933-34, J.~Leray \cite{Ler33, Ler34}
proposed a mathematical question to look for blow-up in
(\ref{NS1}) for $N=3$ driven by the self-similar solutions of the
standard dimensional type, with blow-up time\footnote{Actually,
Leray suggested also to look for a self-similar extension of the
solutions beyond blow-up, i.e., for $t>T$ using the same scaling
similarity variables, \cite[p.~245]{Ler34}.} $T=1$,
 \be
 \label{NS2}
  \tex{
  \vv(x,t)= \frac 1{\sqrt{1-t}} \, {\bf w}(y), \quad p(x,t)= \frac
  1{1-t} \, P(y) \whereA y=\frac x{\sqrt{1-t}}.
 }
 \ee
  Substituting (\ref{NS2}) into (\ref{NS1}) yields for functions ${\bf w}$ and $P$
a 
  ``stationary" system,
 \be
  \label{NS1ss}
   \tex{
\frac 12 \, {\bf w}+  \frac 12 \, (y \cdot \n) {\bf w}+ ({\bf w}
\cdot \n){\bf w}=- \n P + \D {\bf w}, \quad
  \di {\bf w} =0 \inB \ren.
  }
  \ee

During last fifteen years, a number of negative answers
 concerning existence of such non-trivial similarity patterns (\ref{NS2}),
  (\ref{NS1ss}) in $\re^3$ were
 obtained; see \cite{Chae07, Nec96, Mil01},
 and the most advanced and justifying negative answer in \cite{Hou07}.
  Let us note an
  existence result in \cite{Dong07N4} for $N=4$.

\ssk

Nonexistence of Leray's similarity solutions (\ref{NS2}) and other
local types of self-similar blow-up
 is a definite step towards better understanding of
the singularity nature for the Navier--Stokes equations.
Of course,
 this does not settle the problem of singularity formation
 (or nonexistence of finite energy singularities which is more plausible), since there might be other
 ways for (\ref{NS1}) to create singularities as $t \to 1^-$ rather than the purely self-similar
 scenario
  (\ref{NS2}). This multiplicity  question is discussed below.

\ssk

\noi\underline{\em On a countable set of blow-up patterns with
infinite energy}.
 For infinite $L^2$-energy, the  blow-up (also called the enstrophy blow-up of vorticity) in the
 Navier--Stokes equations (\ref{NS1}) can occur even for $N=2$.
 Such {\em global blow-up}  described by {\em von
 K\'arm\'an solutions} is explained in \cite[Ch.~8]{AMGV}, where a rigorous theory
  of such singularities was developed.
  Earlier history of such solutions can be found in \cite{AKPR}.
    This
 blow-up creates a plane jet. It is worth noting that the blow-up
 behaviour {\em is also not} of self-similar form as $t \to T^-$,
 and it  is given by a similarity solutions of a {\em non-local first-order
 Hamilton--Jacobi equation} associated with such a flow. Moreover,
 and this is also crucial, that \cite[pp.~232--235]{AMGV}
  \be
  \label{co1}
   \mbox{there exists a countable set of such different blow-up
   patterns} \quad (N=2).
   \ee
Similar solutions can be constructed for axisymmetric flows in
cylindrical coordinates for the Navier--Stokes equations
(\ref{NS1}) in $\re^3$; see an example in \cite[Ch.~7, \S~3]{AKPR}
and \cite{Oh07}. However, the mathematics of such blow-up patterns
becomes more involved, and (\ref{co1}) for $N=3$ demands difficult
proofs \cite{GalJMP}.

\ssk

\noi\underline{\em Singular (blow-up) set has zero measure}. There
exists another classic direction
 of the singularity theory for the Navier--Stokes equations that
 was originated by Leray himself \cite{Ler34} (details are  available in \cite{Esc03})
  and  in Caffarelli, Kohn, and Nirenberg
 \cite{Caff82}.
 It was shown that the one-dimensional Hausdorff measure of
 the singular (blow-up) points in a time-space cylinder is equal
 to zero. We refer to \cite{Nes07, Ser07} for further development
 and references. In  particular, among other results including Leray's one in \cite{Ler34},
 a refined criterion   is obtained in \cite{Ser07}, saying
 that, if $t=1$ is the first singular (blow-up) moment for a solution
 $\vv(x,t)$ of (\ref{NS1}), then
  \be
  \label{ll1}
   \lim\limits_{t \to 1^-} \frac 1{1-t} \, \displaystyle{ \int\limits_t^1 \int\limits_{\re^3} |\vv(x,t)|^3 \,
   {\mathrm d}x} \, {\mathrm d}t= + \infty.
  \ee
  Performing the scaling as in (\ref{NS2}),
   \be
   \label{ll2}
   \tex{
   \vv(x,t)= \frac 1{\sqrt{1-t}}\, \ww(y,\t), \quad y = \frac x{\sqrt{1-t}}, \quad \t=
   -\ln(1-t) \to +\iy \,\,\, \mbox{as} \,\,\, t \to 1^-,
 }
 \ee
 (\ref{ll1}) takes the same form
\be
  \label{ll3}
   \tex{
   \lim\limits_{\t \to +\infty} \,  {\mathrm e}^{\t}
   \displaystyle{  \int\limits_{\t}^{+\iy} {\mathrm e}^{-s} \Big(
     \int\limits_{\re^3} |\ww(x,s)|^3 \,
   {\mathrm d}y \Big) \, {\mathrm d}s}= + \infty.
 }
  \ee
This means that, for the existence of a singular point  $t=1$, the
solution of the rescaled equations,
 \be
  \label{NS1ssN}
   \tex{
    \ww_\t+
\frac 12 \, {\bf w}+  \frac 12 \, (y \cdot \n) {\bf w}+ ({\bf w}
\cdot \n){\bf w}=- \n P + \D {\bf w}, \quad
  \di {\bf w} =0 \inB \ren,
  }
  \ee
 must diverge (blow-up) as $\t \to +\infty$ in $L^3(\re^3)$.

 Thus, according to the criterion (\ref{ll3}), $t=1$ is not a singular (and hence
 regular) point,
   if the corresponding locally smooth solution of (\ref{NS1ssN})
 does
  not blow-up as $\t \to \infty$ in a suitable functional
 setting. Hence, the problem of global existence and uniqueness
 of a smooth solutions of the Navier--Stokes equations in $\re^3$
 reduces to nonexistence of blow-up in infinite time for the
 rescaled system (\ref{NS1ssN}). In such a framework, this problem falls into the scope of
 standard blow-up/non-blow-up theory for nonlinear evolution PDEs.

  \ssk

  \noi\underline{\em Countable sets of blow-up patterns in
  combustion  problems}.
 Actually, there are many examples of countable sets of blow-up patterns for of much simpler reaction-diffusion
 equations. Amongst them is the classic {\em
Frank-Kamenetskii equation} (1938) \cite{FrK}  developed  in
combustion theory  of solid fuels (also called {\em solid fuel
model}),
  \be
  \label{FK1}
u_t= \D u + {\mathrm e}^u \inA \quad (N=1, \, 2),
 \ee
 for which there exists
 a countable set of different blow-up patterns. Rigorous
 mathematical theory of
 such blow-up patterns was known since the
beginning of the 1990s (see e.g.,  \cite{Vel, GHPV}) and was
  developed by linearization in the inner blow-up region  and
 nonlinear matching. A similar strategy to construct a countable
 set of blow-up patterns is applicable to the higher-order
 reaction-diffusion PDEs \cite{Gal2m}
 \be
  \label{tt3}
   u_t= -(-\D )^m u + |u|^{p-1}u \quad (m \ge 2, \,\,p>1).
    \ee
 where the analysis uses polynomial eigenfunctions and discrete
 spectrum of some related linear non self-adjoint differential operators \cite{Eg4}.

On the other hand, the quasilinear 1D counterpart of (\ref{FK1})
 with the $p$-Laplacian,
 \be
 \label{pp1}
  u_t = (|u_x|^\s u_x)_x+ {\mathrm e}^u \quad (\s>0),
   \ee
is also  known  to admit  a countable set of blow-up patterns
\cite{BuGa} (see also a general discussion in
\cite[pp.~30-34]{GalGeom}), but now, depending on $\s>0$, first
few are self-similar, i.e., represent the case of {\em nonlinear
eigenfunctions}  (not linearized as above for (\ref{FK1})).

\ssk

\noi\underline{\em Evolution completeness as a necessary
ingredient}. It is a principal open problem to describe the whole
set of all possible blow-up patterns (if any) for the
Navier--Stokes equations (\ref{NS1}). Evidently, proving
nonexistence of all the blowing up patterns, i.e., nonexistence of
blow-up at all, will settle
 the fundamental problem of global smooth (and hence unique)
 continuation of sufficiently arbitrary solutions.
Another natural possibility is to establish that, for given class
of data, the orbits do not approach any of blow-up pattern
scenario, so remain regular for all times.

\ssk

In this context, the problem of {\em evolution completeness} of
the given countable set of patterns occur, meaning that these
patterns exhaust all possible types of approaching the singularity
for a fixed class of data. For linear problems,
  the evolution completeness  follows from
standard completeness and closure of eigenfunction subset of a
linear operator or pencil in a fixed functional framework, so the
evolution completeness does not have a separate meaning. For
nonlinear problems, where ``linear" notions of completeness and
closure make no sense in general, the evolution completeness
becomes key.

\ssk

There are some examples of evolution completeness of countable
sets of linearized or nonlinear patterns (eigenfunctions) in
parabolic asymptotic theory. For instance, the full classification
of blow-up sets for the Frank-Kamenetskii equation (\ref{FK1}) was
performed by Vel\'azquez \cite{Vel} by actual proving the
evolution completeness of the countable set of linearized blow-up
patterns. Concerning nonlinear patterns, it seems that there
exists a unique example of a proof of evolution completeness of
such a countable set for the porous medium and $p$-Laplacian
equations in 1D or in radial geometry in $\ren$,
 $$
 u_t = \D (|u|^{m-1} u) \andA u_t= \n \cdot (|\n u|^{m-1} \n u)
 \withA m>1;
 $$
 see \cite{CompG},
 where the notion of evolution completeness was
introduced.

\ssk

\noi\underline{\em On the well-posed Burnett-type
equations.}
 Concerning our conclusion, as a straightforward application of the
scaling technique in {\bf (IV)},
we present a simple proof of the fact that a local smooth solution
$\vv(x,t)$ of (\ref{NS1}), which is {\em uniformly bounded in
$L^p(\ren)$}, with
  \be
  \label{NS3}
 p > N, \quad \mbox{i.e.,  \, $p > 3$ \, for \, $N=3$},
   \ee
  is uniformly bounded in
  $L^\infty(\ren)$,  so that such $L^p$-solutions do not blow-up.
  The non-blow-up also is proved in the well-known critical case $p=N=2$.
  The condition (\ref{NS3}) is consistent with
Leray--Prodi--Serrin--Ladyzhenskaya regularity  $L^{p,q}$ criteria
 and other more recent
  results; see
  key references, history, details, and results in recent papers \cite{Qio07, Esc03, Ser07, Gala07}.
 In addition, (\ref{NS3}) re-derives  Kato's class of unique mild
 solutions in $\ren$,
 \cite{Kato84}; see details and key references in \cite{Gala07, Way05}.

  This approach also extends to the related $2m$th-order {\em well-posed Burnett
  equations},
 \be
  \label{NS1m}
 \vv_t +(\vv \cdot \n)\vv=- \n p - (-\D)^m \vv, \quad \di \vv=0
 \inA,
  \ee
 containing the
    higher-order diffusion operator $-(-\D)^m \vv$, with arbitrary $m \ge
    1$.
   Then the regularity criterion similar to that in (\ref{NS3}) for $m=1$
   reads
   $$
   \tex{
   p > \frac N{2m-1} \forA m=1,2,3,... \, .
  }
  $$

On the other hand, we prove that for smooth fast decaying
divergence-free $L^2$-data $\vv_0$, finite-time blow-up is
impossible in dimensions
 $$
 N < 2(2m-1) \quad (\mbox{$N =2$ for $m=1$ is included, \cite{Ler33N2, Lad70}}),
 $$
 so there exists a unique global classic bounded solution.



\section{Method of interpolation: global
existence for the KSE}
 \label{SGl}

\subsection{A priori estimates}

To demonstrate these classic approaches, we  consider  the one dimensional KSE
in the  following IBVP setting (for convenience, here $D= \frac
{\partial}{{\partial} x}$):
\begin{equation}
\label{e2.1}
 \tex{
v_t+D^4 v+D^2 v=\frac{1}{2}\, D v^2, \quad t>0, \quad x\in(-L,L),
 }
\end{equation}
\begin{equation}\label{e2.2}
v(x,0)=v_0(x), \quad x\in(-L,L),
\end{equation}
with either:
\begin{equation}\label{e2.3}
v=Dv=0, \quad x=-L, \quad x=L, \quad t>0 \orA
\end{equation}
\begin{equation}\label{e2.3'}
v=D^2v=0, \quad x=-L, \quad x=L, \quad t>0.
\end{equation}

\ssk

 \noi\underline{\em A priori estimates}.
 Multiplying (\ref{e2.1}) by $v$ and
integrating by parts over $\Omega=(-L,L)$ with regard to
(\ref{e2.3}) or (\ref{e2.3'}), we find
\begin{equation}
\label{e2.4}
\tex{
 \frac 12 \,\frac{{\mathrm d}}{{\mathrm d}t}\displaystyle{\int\limits_{\Omega}
v^2(x,t)\,{\mathrm d}x+\int\limits_{\Omega} |D^2 v|^2\,{\mathrm
d}x-\int\limits_{\Omega} |Dv|^2\,{\mathrm d}x=0.}
 }
\end{equation}
Due to (\ref{e2.3}) or (\ref{e2.3'}), we have
\begin{equation}\label{e2.5}
 \tex{\displaystyle{\int\limits_{\Omega} |Dv|^2\,{\mathrm d}x} \le \Big(\displaystyle{\int\limits_\O
v^2\,{\mathrm d}x\Big)^{\frac 12} \Big(\int\limits_{\Omega} |D^2
v|^2\,{\mathrm d}x\Big)^{\frac 12}. }}
\end{equation}
Denoting by
\[
 \tex{
E(t)=\displaystyle{\int\limits_{\Omega} v^2(x,t)\,{\mathrm d}x},
 }
\]
from (\ref{e2.4}) we get
\[
 \tex{
\frac 12 \,\displaystyle{\frac{{\mathrm d}E}{{\mathrm d}t}} \le
\Big(\int\limits_{\Omega} |D^2 v|^2\,{\mathrm d}x\Big)^{\frac
12}E^{\frac 12}(t)-\int\limits_{\Omega} |D^2 v|^2\,{\mathrm d}x
\le \tex{\frac 14} \,E(t).
  }
\]
So
\begin{equation}\label{e2.6}
E(t)\le {\hat E}(t):=E(0){\mathrm e}^{ \frac t2}.
\end{equation}
Next, integrating (\ref{e2.4}) in $t>0$, we see that
\begin{equation}\label{e2.7}
 \tex{
\displaystyle{\int\limits_0^t \int\limits_{\Omega} |D^2
v|^2\,{\mathrm d}x}\,{\mathrm d}t = \int\limits_0^t
\int\limits_{\Omega} |Dv|^2\,{\mathrm d}x\,{\mathrm d}t -
\tex{\frac 12 \,E(t)+\frac 12 \,E(0)}.
 }
\end{equation}
Due to (\ref{e2.3}) or (\ref{e2.3'}) we have
\[
 \tex{
0=\displaystyle{\int\limits_0^t \int\limits_{\Omega} D(v Dv)\,{\mathrm
d}x\,{\mathrm d}t = \int\limits_0^t \int\limits_{\Omega} v D^2
v\,{\mathrm d}x\,{\mathrm d}t +\int\limits_0^t
\int\limits_{\Omega} |Dv|^2\,{\mathrm d}x\,{\mathrm d}t.
 }}
\]
From this, it follows that
\begin{equation}
  \label{e2.8}
 \int_0^t \int_{\Omega}
|Dv|^2\,{\mathrm d}x\,{\mathrm d}t \le \frac 12 \, \int_0^t
\int_{\Omega} v^2\,{\mathrm d}x\,{\mathrm d}t + \frac 12 \,
\int_0^t \int_{\Omega} |D^2 v|^2\,{\mathrm d}x\,{\mathrm d}t.
\end{equation}
Thus, (\ref{e2.7}) implies that
\[
\frac 12 \,\displaystyle{ \int\limits_0^t \int\limits_{\Omega}
|D^2 v|^2\,{\mathrm d}x\,{\mathrm d}t \le \frac 12 \,
\int\limits_0^t E(t)\,{\mathrm d}t- \frac 12 \,E(t)+ \frac 12
\,E(0). }
\]
Hence, by (\ref{e2.6}) one has
\begin{equation}\label{e2.9}
 \tex{
\displaystyle{\int\limits_0^t \int\limits_{\Omega} |D^2 v|^2\,{\mathrm
d}x\,{\mathrm d}t \le 2E(0){\mathrm e}^{\frac t2}-E(0).
 }}
\end{equation}
Then from (\ref{e2.7}) it follows that
\begin{equation}\label{e2.10}
 \tex{
\displaystyle{\int\limits_0^t \int\limits_{\Omega}
|Dv|^2\,{\mathrm d}x\,{\mathrm d}t \le \big(2{\mathrm e}^{\frac
t2}-\tex{\frac 32}\big)E(0).
 }}
\end{equation}

\ssk

\noi\underline{\em  Estimate of} $\displaystyle{\int\limits_0^t
{\int\limits_{\Omega} v^2 |Dv|^2\,{\mathrm d}x\,{\mathrm d}t}}$.
Thanks to embedding inequality for $\Omega \subset \re$ we have
\[
 \tex{
\|v(t)\|^2_{\infty} \le c_{\infty}
\displaystyle{\int\limits_{\Omega}|Dv|^2(x,t)\,{\mathrm d}x
 }}
\]
since $v|_{\prt \Omega}=0$. Then
\[
\tex{ \displaystyle{ \int\limits_{\Omega} v^2 |Dv|^2\,{\mathrm
d}x} \le c_{\infty}
\displaystyle{\Big(\int\limits_{\Omega}|Dv|^2\,{\mathrm
d}x\Big)^2.
 }}\]
Next
\[
\tex{\displaystyle{\Big(\int\limits_{\Omega}|Dv|^2\,{\mathrm
d}x\Big)^2 \le \int\limits_{\Omega}v^2\,{\mathrm d}x \cdot
\int\limits_{\Omega}|D^2 v|^2\,{\mathrm d}x = E(t)
\int\limits_{\Omega}|D^2 v|^2\,{\mathrm d}x.
 }}
\]
Applying (\ref{e2.9}), we find
\begin{equation}\label{e2.11}
\begin{array}{ll}
\mbox{
 $
\displaystyle\int\limits_0^t \int\limits_{\Omega} v^2|Dv|^2\,{\mathrm d}x\,{\mathrm d}t
\le c_{\infty} \int\limits_0^t {\hat E}(\tau)
\int\limits_{\Omega}|D^2 v|^2\,{\mathrm d}x\,{\mathrm d}\tau
 $
 }
\\
\tex{
 \le c_{\infty} {\hat E}(t) \displaystyle{\int\limits_0^t}
\int\limits_{\Omega}|D^2 v|^2\,{\mathrm d}x\,{\mathrm d}\tau \le
c_{\infty} {\hat E}(t)\big(2{\mathrm e}^{\frac t2}-1\big)E(0)
 } \ssk\ssk
\\
 \tex{
\dip{=c_{\infty}{\mathrm e}^{\frac t2}\big(2{\mathrm e}^{\frac
t2}-1\big)E^2(0).}
 }
\end{array}
\end{equation}
Next, by using the above estimates and the linear theory of the
parabolic equations of the 4th-order we get (here
$Q_T=\Omega\times(0,T)$)
\begin{equation}\label{e2.12}
\begin{array}{ll}
 \tex{
\dip{ \|v(T)\|^2_{L^2(\Omega)}+\|v_t\|^2_{L^2(Q_T)} + \|D^4
v\|^2_{L^2(Q_T)}}
 }
 \ssk\ssk
\\
 \tex{
\dip{ +\|D^2 v\|^2_{L^2(Q_T)}+\|Dv\|^2_{L^2(Q_T)} +
\|v\|^2_{L^2(Q_T)}}
 } \ssk\ssk
\\
 \tex{
 \dip{\le
C(\|vDv\|^2_{L^2(Q_T)}+\|v_0\|^2_{L^2(\Omega)})}
 }\ssk\ssk
\\
\tex{ \dip{\le C_1\big({\mathrm e}^TE^2(0)+{\mathrm e}^{\frac T
2}E(0)+\|v_0\|^2_{L^2(\Omega)}\big)}
 }
\end{array}
\end{equation}
with a positive constant $C_1$ independent of $v$ and $T>0$.


\subsection{Global existence}

\begin{theorem}
\label{Th.3.1}
For any $v_0\in L^2(\Omega)$, both the IBVPs
$(\ref{e2.1})$--$(\ref{e2.3})$ and $(\ref{e2.1})$, $(\ref{e2.2})$,
$(\ref{e2.3'})$ have solutions for any $t>0$, satisfying
inequality $(\ref{e2.12})$.
 \end{theorem}


\noi{\em Proof.} The ground {\em a priori} estimates (\ref{e2.6}),
(\ref{e2.8}), (\ref{e2.9}) and (\ref{e2.11}) are obtained by
multiplication of equation (\ref{e2.1}) by $v$. Thanks to this we
use the Galerkin approach and these {\em a priori} estimates are
valid for the Galerkin approximations $\{v_m\}$.

The passage to the limit $v_m\to v$ as $m\to\infty$ in the
nonlinear term $v(x,t)Dv(x,t)\ (x\in\Omega\subset\re,\,t>0)$ is
proved by standard arguments (see \cite{1, JLi}). As a result we
get a weak solution. Then from (\ref{e2.11}) due to the linear
theory of parabolic equations of the 4th order we obtain the final
estimate (\ref{e2.12}). $\qed$

\section{Method of eigenfunctions: blow-up of solutions}
 \label{SBlowup}

We now show how, using a similar interpolation-eigenfunction
technique, to derive sufficient conditions of blow-up for the KSE
 with special ``boundary conditions".

\subsection{Basic computations}

 Multiplying equation (\ref{e2.1})
by a function $\psi(x)$ belonging to $C^4(\re)$, we get
\begin{equation}\label{e4.3}
 \tex{
\frac{\mathrm d}{{\mathrm d}t} }
\displaystyle{\int\limits_{\Omega}v\psi\,{\mathrm d}x
+\int\limits_{\Omega}v(D^4\psi+D^2\psi)\,{\mathrm d}x
=-\mbox{$\frac 12$} \,\int\limits_{\Omega}v^2D\psi\,{\mathrm
d}x+B(v,\psi), }
\end{equation}
where $\Omega=(0,L)\subset\re$ and
\begin{equation}\label{e4.4}
 \tex{
B(v,\psi)=\big[vD^3\psi-Dv\cdot D^2\psi+(D^2 v+v)D\psi-(D^3
v+Dv)\psi+\frac 12 \,v^2\psi\big]^L_0. }
\end{equation}
For the function
\[
\psi(x)=\psi_{\lambda}(x):=|x-L|^{\lambda}
\]
with $\lambda>6$, we have from (\ref{e4.3})
\begin{equation}\label{e4.5}
 \tex{
\frac{\mathrm d}{{\mathrm d}t}}
\displaystyle{\int\limits_{\Omega}v\psi_{\lambda}\,{\mathrm d}x =
\mbox{$\frac{\lambda}{2}$}\, \int\limits_{\Omega}v^2\cdot
|x-L|^{\lambda-1}\,{\mathrm
d}x-\int\limits_{\Omega}v\psi_1\,{\mathrm d}x}+B_0(v) \whereA
\end{equation}
\[
 \tex{
\psi_1=\lambda(\lambda-1)|x-L|^{\lambda-4}[(\lambda-2)(\lambda-3)+|x-L|^2]
\quad \mbox{and} }
\]
\begin{equation}\label{e4.6}
\begin{array}{ll}
\dip{
 \tex{
B_0(v)=-\frac 12 \,L^{\lambda}v^2(0,t)+\lambda
L^{\lambda-1}\cdot[(\lambda-1)(\lambda-2)L^{-2}]v(0,t) }}
 \ssk\ssk\ssk\ssk \\
\dip{+L^{\lambda}[1+\lambda(\lambda-1)L^{-2}]Dv(0,t)+ \lambda
L^{\lambda-1}D^2 v(0,t)+L^{\lambda}D^3 v(0,t).}
\end{array}
\end{equation}
Denote by
\[
 \tex{
J(t)=\displaystyle{\int\limits_{\Omega}v(x,t)|x-L|^{\lambda}\,{\mathrm d}x},
\quad \Omega=(0,L). }
\]
Thanks to the inequality
\[
\begin{array}{ll}
\dip{
\displaystyle\int\limits_{\Omega}v\psi_1\,{\mathrm d}x}
 &
\dip{\le  \Big(\int\limits_{\Omega} v^2
|x-L|^{\lambda-1}\,{\mathrm d}x\Big)^{\frac 12}
\Big(\int\limits_{\Omega}\frac{\psi_1^2(x)}
{|x-L|^{\lambda-1}}\,{\mathrm d}x\Big)^{\frac 12}}
 \ssk\ssk
 \\
 & \dip{\le \mbox{$\frac{\lambda}{4}$}
\displaystyle{ \int\limits_{\Omega} v^2
 |x-L|^{\lambda-1}\,{\mathrm d}x}+C_{\lambda}(L)} \withA
\end{array}
\]
\[
 \tex{
C_{\lambda}(L)=\lambda(\lambda-1)L^{\lambda-6}\big[
\frac{(\lambda-2)^2 (\lambda-3)^2}{\lambda-6}+ \frac{2(\lambda-2)
(\lambda-3)}{\lambda-4}L^2+\frac{1}{\lambda-2}L^4\big], }
\]
from (\ref{e4.5}) it follows that
\[
\tex{ \displaystyle\frac{{\mathrm d}J}{{\mathrm d}t}\le
\tex{\frac{\lambda}{4}}\, \displaystyle\int\limits_{\Omega}
v^2|x-L|^{\lambda-1}\,{\mathrm d}x+B_0(v)-C_{\lambda}(L).
 }
\]
By using the H\"older inequality
\[
 \tex{
J^2(t)\le \frac{L^{\lambda+2}}{\lambda+2}\displaystyle\int\limits_{\Omega}
v^2|x-L|^{\lambda-1}\,{\mathrm d}x,
 }
\]
we get
\begin{equation}\label{e4.7}
\tex{\displaystyle
 \frac{{\mathrm d}J}{{\mathrm d}t}\ge
\tex{\frac{\lambda(\lambda+2)}{4}}\,
\frac{1}{L^{\lambda+2}}J^2(t)+B_0(v)-C_{\lambda}(L).
 }\end{equation}

\subsection{General initial boundary value problem}

We consider a non-standard IBVP for the KSE (\ref{e2.1}) that
includes
 initial data
\begin{equation}\label{e5.2}
v(x,0)=v_0(x)
\end{equation}
and  some general boundary conditions to be specified,
\begin{equation}\label{e5.3}
B_0(v)=h(t) \quad \mbox{at $x=0$ for $t>0$.}
\end{equation}

We understand a solution $v$ of (\ref{e2.1}), (\ref{e5.2}),
(\ref{e5.3}) as a weak solution in the sense of identity
(\ref{e5.3}) with respect to any test function $\psi\in C^4(\re)$
with an initial data in the sense
\[
 \tex{\displaystyle
\int\limits_{\Omega} v(x,t)\psi(x)\,{\mathrm d}x \to
\int\limits_{\Omega} v_0(x)\psi(x)\,{\mathrm d}x } \asA t \to 0.
\]
Of course, the function $v$ is assumed to belong to the function
space such that identity (\ref{e5.3}) makes sense.

\subsection{The blow-up results}

As follows from (\ref{e4.7}) the blow-up phenomenon for IBVP
(\ref{e2.1}), (\ref{e5.2}), (\ref{e5.3}) depends on
\begin{equation}\label{e6.1}
H_{\lambda}(v,L):=B_0(v)-C_{\lambda}(L).
\end{equation}
It means it depends on the relationship between initial data
$v_0(x)$ and boundary conditions.


 \begin{theorem}
  \label{Th.6.1}
 Let $H_{\lambda}(v,L)\ge a^2>0$ for some
$\lambda>6$ and $L>0$. Then there is no global  solution of
$(\ref{e2.1})$, $(\ref{e5.2})$, $(\ref{e5.3})$.

 Moreover, there is no global
 solution of $(\ref{e2.1})$, $(\ref{e5.2})$, $(\ref{e5.3})$ for
$x\in(0,L)$, and
\[
 \tex{
J(t)\ge \frac{a}{\kappa}\tan(a\kappa t+c_0) \whereA }
\]
\[
 \tex{
 \kappa=\big(\frac{\lambda(\lambda+2)}{4}L^{-(2+\lambda)}
\big)^{\frac 12},} \,\,\, a=\sqrt{a^2}>0,\,\,\,
 \tex{
c_0=\arctan\frac{\kappa J_0}{a},}\,\,\,  {\dip \tex{
J_0=\displaystyle\int\limits_{\Omega}v_0(x)|x-L|^{\lambda}\,{\mathrm d}x,}}
\]
and  the blow-up time $T_{\infty}$ is estimates as:
\[
 \tex{
T_{\infty}\le \frac{\frac \pi 2-c_0}{a\kappa}.
 }
\]
 \end{theorem}

The proof follows immediately from (\ref{e4.7}).


\begin{corollary}
 \label{Cor.6.1}
Suppose that  $v(x,t)$  satisfies,
\[
\begin{array}{ll}
v(0,t)=0,\\
   Dv(0,t)+D^3v(0,t)\ge {\rm const}>0 \quad \forall
t>0, \,\,\mbox{and}\\
 D v(0,t), \,\, D^2 v(0,t) \ge 0.
\end{array}
\]
Then there is $L_0>0$ such that all the assumptions
of Theorem $\ref{Th.3.1}$  are  fulfilled.
 \end{corollary}


\begin{theorem}
 \label{Th.6.2}
 Let $H_{\lambda}(v,L)\ge 0$ for some
$\lambda>6$ and $L>0$. Let
\[
 \tex{
J_0:=\displaystyle\int\limits_{\Omega}v_0(x)|x-L|^{\lambda}\,{\mathrm d}x>0.
 }
\]
Then there is no global solution of $(\ref{e2.1})$,
$(\ref{e5.2})$, $(\ref{e5.3})$.

Moreover, there is no global  solution of $(\ref{e2.1})$,
$(\ref{e5.2})$, $(\ref{e5.3})$ for $x\in(0,L)$. In particular $J$ satisfies
\[
 \tex{
  J(t)\ge \frac{J_0}{1-J_0 \kappa^2 t},
 }
\]
so that  for the blow-up time $T_{\infty}$ we have
\[
 \tex{
T_{\infty}<\frac{1}{\kappa^2 J_0}.
 }
\]
\end{theorem}

The proof follows  from (\ref{e4.7}).

\medskip

\begin{theorem}
 \label{Th.6.3}
Let $H_{\lambda}(v,L)\ge -a^2$ for some
$\lambda>6$ and $L>0$. Let
\[
 \tex{
J_0:=\displaystyle\int\limits_{\Omega}v_0(x)|x-L|^{\lambda}\,{\mathrm
d}x>\frac{|a|}{\kappa}.
 }
\]
Then there is no global solution of $(\ref{e2.1})$,
$(\ref{e5.2})$, $(\ref{e5.3})$.

Moreover, there is no global solution of $(\ref{e2.1})$,
$(\ref{e5.2})$, $(\ref{e5.3})$ for $x\in(0,L)$.  In particular
\[
 \tex{
J(t)\ge \frac{1+c_0 {\mathrm e}^{2a\kappa t}}{1-c_0 {\mathrm
e}^{2a\kappa t}}\frac{a}{k} \withA
 }
 \tex{
c_0=\frac{\kappa J_0-a}{\kappa J_0+a}.
 }
\]
Hence, the blow-up time $T_{\infty}$ satisfies:
\[
 \tex{
T_{\infty}<\frac{1}{2a\kappa}\ln \frac{\kappa J_0+a}{\kappa J_0-a}
\whereA
 }
 \tex{
\kappa=\big(\frac{\lambda(\lambda+2)}{4}L^{-(2+\lambda)}\big)^{\frac
12}
 } \quad \mbox{and $a=|a|$}.
\]
 \end{theorem}

Again, the proof  follows directly  from (\ref{e4.7}).

\medskip

\begin{corollary}
 \label{Cor.6.3}
 Let $v_0(x)\ge c|x|^{\mu}$ with some
$\mu>0$ and $c>0$. Suppose that the boundary values $v,\,Dv,\,D^2v$, and
$D^3 v$ at $x=0$ do not depend on $t>0$. Then there exists $\lambda>0$
such that the assumptions of Theorem $\ref{Th.6.3}$ are fulfilled.
 \end{corollary}

\section{Global existence and $L^\iy$-bounds by  weighted Gronwall's
inequalities}
 \label{SGlH}

\subsection{A general KS-type model}

In this section, we extend  global existence approaches to more
general $2m$th-order parabolic equations of the KS-type.
 As it has been seen, the zero boundary conditions of the IBVPs
 under consideration reinforced the energy control and helped the global solvability of the KSE.
It is then  reasonable  to consider the KS-type equations in
unbounded domains suggesting  ``infinite propagation" with no
spatial ``obstacles" and bounds, which the Cauchy problem is most
natural for.

 Thus, 
  as a basic and typical
 model,
we consider the Cauchy problem (CP) for the equation, which
includes both stable and unstable linear diffusion terms as well
as the ``convection":
 \be
 \label{h1}
 v_t=-(-\D)^{2l} v + (-\D)^l v+  \BB_1 (|v|^p) \inA  \quad (m=2l= 2,4,6,...),
  \ee
  where $p > 1$ and  $\BB_1$ is the linear first-order differential
  operators in divergence form
   \be
   \label{h2}
   \tex{
   \BB_1 u= \frac 1{p} \, \sum_{(k)} d_k D_{x_k} u.
  }
  \ee
 For simplicity,  the coefficients $\{d_k\}$
of the operator are assumed to be constant, though some
$(x,t)$-dependence can be easily allowed.
 The classic KSE (\ref{e2.1}) corresponds to
 \be
 \label{KSpar}
 p=2, \quad l=1, \andA N=1 \quad (d_1=1).
  \ee
  We refer to (\ref{h1}) as to the {\em modified Kuramoto--Sivashinsky
 equation}
 (the mKSE).

Thus, we consider for (\ref{h1}) the CP with bounded sufficiently
smooth  data
 \be
 \label{h3}
  v(x,0)=v_0(x) \inB \ren, \quad v_0 \in L^\infty(\ren) \cap
  H^{2m}(\ren).
   \ee
For instance, for the main demonstrations of the techniques
involved, we may always assume that $v_0(x)$ is smooth and has
exponential decay at infinity.

One of the key mathematical features of the mKSE (\ref{h1}) is
that it admits multiplication by $v$ in $L^2$ to get the {\em a
priori} $L^2$-bound. Namely, similar to related straightforward
manipulations in Section \ref{SGl}, we have by the H\"older
inequality:
 \be
 \label{h3s}
   \begin{matrix}
   \frac 12\,
  \frac{{\mathrm d}}{{\mathrm d}t} \|v(t)\|_2^2=
  - \|\D^l v(t)\|_2^2 +\displaystyle \int v(-\D)^l v \,\,dx\le
 - \|\D^l v(t)\|_2^2 + \|v(t)\|_2 \|\D^l v(t)\|_2\qquad
  \ssk\ssk\ssk\\
  \le \frac 14\,
 \|v(t)\|_2^2 \quad
  \Longrightarrow \quad
  \|v(t)\|_2^2 \le \|v_0\|_2^2 \, {\mathrm e}^{\frac t2} \quad \mbox{for all} \quad t>0.
  \qquad
 \end{matrix}
   \ee

Another obvious higher-order model with a similar $L^2$-control
(\ref{h3s}) is
 \be
 \label{h1aa}
  \tex{
 v_t=-(-\D)^{m} v + \frac 14\,  v+  \BB_1 (|v|^p) \inA  \quad (m= 2,3,4,...),
  }
  \ee
with the same convection-nonlinear dispersion and the linear
unstable zero-order term $+ \frac 14\, v$.
 On the other hand, $\frac 14\, v$ can be replaced by $-\D v$ as in (\ref{e2.100}).
  However, for some
$2m$th-order KS-type models, deriving $L^2$-bounds can represent a
hard problem, so we avoid using those in future application of the
scaling and other techniques.

  Thus, as a principal issue of our analysis,
 the  intention is to use this {\em a priori} estimate (\ref{h3s}) to
 derive stronger uniform bounds on solutions $v(x,t)$ in the $\sup_x$-norm
 for all $t \ge 0$.
 The developed  scaling technique is rather general and applies to
other higher-order KS-type models with {\em a priori} known $L^2$,
 $L^p$, Sobolev, or other types of weaker or stronger integral bounds on solutions.

\subsection{Fundamental solution}

We will need the fundamental solution of the corresponding linear
parabolic (poly-harmonic) equation
 \be
  \label{Lineq}
  u_t = - (-\D)^m u \quad \mbox{in} \,\,\, \ren \times \re_+ \quad
  (m=2l),
  \ee
which has the self-similar form
  \be
  \label{1.3R}
   b(x,t) = t^{-\frac N{2m}}F(y), \quad y= x/t^{\frac 1{2m}}.
 \ee
 The rescaled kernel $F$ is the unique radial solution of the elliptic equation
  \begin{equation}
\label{ODEf}
 \mbox{$
 {\bf B} F \equiv -(-\Delta )^m F + \frac 1{2m}\, y \cdot
\nabla F  + \frac N{2m}\, F  = 0 \quad \mbox{in} \,\,\, \ren,
 \quad \int F = 1.
 $}
\end{equation}
$F$ is oscillatory for $m >1$, has exponential decay as $|y| \to
\infty$, and satisfies, for some positive constants $D$ and $d$
depending on $m$ and $N$
 \cite{EidSys},
\begin{equation}
\label{fbar}
 \mbox{$
 |F(y)| < D \, {\mathrm e}^{-d|y|^{\alpha}}
\inB \ren \whereA \a=  \frac {2m}{2m-1} \in (1,2).
 $}
\end{equation}

\subsection{Local existence and uniqueness of smooth solutions}

 As a standard practice  (see
 \cite{EidSys, Fr} and \cite[Ch.~15]{Tay}),  for sufficiently regular initial data,
{\em local} in time existence and uniqueness of the classical
solution of the CP (\ref{h1}),
(\ref{h3}) is studied (via Duhamel's principle)  by using  the 
 integral equation, where we have integrated by parts once,
 \be
 \label{1.2N}
  \tex{
  v(t) =
     b(t)*v_0 +
     \displaystyle\int\limits_0^t (-\D)^l b(t-s)* v(s)\, {\mathrm d}s
     + \int\limits_0^t \BB_1^* \, b(t-s)* |v|^{p}(s)\, {\mathrm d}s,
     }
 \ee
 where $b(t)$ is the fundamental solution (\ref{1.3R}) and $\BB_1^*$ is the adjoint operator,
  \be
  \label{B**}
 \tex{
 \BB_1^*(\cdot)
  = - \frac 1p\, \sum_{(k)}D_{x_k}(d_k (\cdot)).
 }
 \ee
For more singular  data $v_0 \in L^q$, with  $q \in (1,\infty)$,
  the solution may  not be classical and is then understood as a proper continuous curve $u:[0,T]
\to L^q$  satisfying equation (\ref{1.2N}). Such questions of
local existence and uniqueness were first systematically studied
by Weissler in the 1970s and 80s, \cite{Weis0, Weis}.
 See the results in \cite[pp.~87-90]{Weis}, which actually  apply to 
  $2m$-th order equations like (\ref{h1}).
  More recent 
 results on local and global existence for higher-order parabolic equations
 such as
 (\ref{1.2N}) can be found in \cite{Bao, Cui, Eg4, GP1}.

In  what follows, we are interested in global existence and
behaviour of solutions, so we have assumed that  data $v_0$ are
sufficiently regular, and then we will use the equivalent integral
equation (\ref{1.2N}) for such purposes.

\subsection{Global existence by Henry's version of weighted
Gronwall's inequality}


\begin{theorem}
 \label{Th.G2}
 The Cauchy problem $(\ref{h1})$, $(\ref{h3})$ has a global unique
 classical solution if
 \be
  \label{h10}
 \tex{
1<  p \le 3 \andA  p < p_0= 1+ \frac{2(2m-1)}{N}.
  }
   \ee
     \end{theorem}

     For the original KSE (\ref{e2.1}) with parameters
     (\ref{KSpar}), (\ref{h10}) is valid in dimension
      \be
      \label{dd1}
      N  < 6.
       \ee

\noi{\em Proof.} Let us write (\ref{1.2N}) in greater detail,
\be
\label{h5}
  \begin{matrix}
  v(x,t)= t^{-\frac N{2m}}\displaystyle \int\limits_{\ren}
  F(\cdot)v_0(y)\, {\mathrm d}y+
  \displaystyle\int\limits_0^t(t-s)^{-\frac {N+2l}{2m}}\, {\mathrm d}s
  \int\limits_{\ren} (-\D)^l F(\cdot)
 v(y,s) \, {\mathrm d}y\qquad
  \ssk\ssk\\
  + \,
   \displaystyle\int\limits_0^t(t-s)^{-\frac {N+1}{2m}}\, {\mathrm d}s \int\limits_{\ren}
    \BB_1^* F\big(\tex{\frac{x-y}{(t-s)^{1/2m}}}\big)
 |v|^p(y,s) \, {\mathrm d}y, \quad
 (\cdot)=\big(\tex{\frac{x-y}{(t-s)^{1/2m}}}\big),\qquad
 \end{matrix}
  \ee
where we mean in $(-\D)^lF(z)$ and  $\BB_1^* F(z)$ that the
operators act in the $z$-variable.
 Denote
 \be
  \label{h6}
  V(t)= \sup_{x \in \ren} \, |v(x,t)|.
   \ee
 Writing in (\ref{h5}) $|v|^p= |v| \, |v|^{p-1}$ for $p \le 3$ ($p-1 \le
 2$),
 using
 H\"older's inequality yields
   \be
 \label{1.3N}
 \begin{matrix}
|v(t)| 
 \le \sup |v_0| \int |b(t)| +
\displaystyle\int\limits_0^t | \D^l b(t-s)
* v(s)| \, {\mathrm d}s
\ssk\ssk\qquad\\
 +\,
  \displaystyle\int\limits_0^t | \BB_1^* b(t-s)
* (|v|^{p-1}|v|)(s)| \, {\mathrm d}s
\le \, \sup |v_0| \, \|F\|_1
  \ssk\ssk\qquad\\
 + \,
 C \|\D^l F\|_{1}
\displaystyle\int\limits_0^t V(s) (t-s)^{-\frac l m} \, {\mathrm
d}s \ssk\ssk\qquad\\
 + \,
 C \|B_1^* F\|_{\frac 2{3-p}}
\displaystyle\int\limits_0^t V(s)\|v(s)\|_2^{p-1} (t-s)^{\b-1} \, {\mathrm
d}s,\qquad
\end{matrix}
  \ee
 \be
 \label{h7}
  \tex{
   \mbox{where} \quad
 \b= \frac{4m-2 -N(p-1)}
 {4m}.
  }
  \ee

 Thus, by (\ref{h3s}), we obtain weighted    Gronwall's inequality
 (here $\frac l m=\frac 12$)
 \be
 \label{h8}
 \mbox{$
 V(t) \le C + C \displaystyle\int\limits _0^t (t-s)^{- \frac 12} V(s) \, {\mathrm d}s
 + C \displaystyle \int\limits _0^t {\mathrm e}^{\frac{(p-1)s}4}(t-s)^{\b-1} V(s) \, {\mathrm
 d}s,
 $}
 \ee
 where, obviously,  the last term on the right-hand side is key for boundedness of $V(s)$.
 By Henry's estimates for such weighted inequalities
\cite[p.~188]{He}, it follows  that $V(s)$ is bounded for any
$t>0$,
if  $\b
> 0$, which is equivalent to the last inequality in (\ref{h10}).
$\qed$


\subsection{On a double exponential $L^\iy$-growth by weighted Gronwall's
inequality}

As we have mentioned, for the KSE, it is principal also to
establish the best estimate on the growth of global solutions for
$t \gg 1$, to be compared with the exponential
one
(\ref{aa1}).

Let us see what kind of $L^\infty$-bound is guaranteed by the
approach applied above. Consider the principal integral operator
in (\ref{h8}), where we skip all the constants $C$,
 \be
 \label{h8n}
 \mbox{$
 V(t) \le 1 + \displaystyle \int\limits _0^t {\mathrm e}^{\frac{(p-1)s}4}(t-s)^{\b-1} V(s) \, {\mathrm
 d}s \quad (\b>0).
 $}
 \ee
For $p=1$, Henry's ``discrete" proof in \cite[p.~188]{He} also
gives the fact that solutions of such Gronwall's inequalities do
not grow as $t \to \infty$ not faster than exponentially, which is
fine according to (\ref{aa1}).

With the presence of the multiplier ${\mathrm e}^{\frac{(p-1)s}4}$
in the kernel in (\ref{h8n}), this is not the case. Using the idea
of Henry's proof, a certain estimate of the behaviour of $V(t)$
for $t \gg 1$ can be indeed obtained. Nevertheless, rather
surprisingly,  this will not be an exponential growth, which may
naturally look most plausible. Namely,  it is easy to see that an
upper bound on the growth of solutions of (\ref{h8n}) is {\em
doubly exponential} in the sense that their solutions cannot grow
as $t \to +\iy$ faster than the ``supersolution" (corresponding to
``$=$" in (\ref{h8n}))
 \be
 \label{ff1}
  \tex{
  \hat V(t) = \exp\big\{ \big[\frac 4{\b(p-1)}+\e \big] \, t^\b\, {\mathrm
  e}^{\frac{(p-1)t}4}\big\}, \quad \mbox{with arbitrarily small} \,\,\,
  \e>0.
   }
   \ee
 In other words, the approach based on weighted Gronwall's
 inequalities, though proving global existence, supplies us with
a  rather non-realistic (in comparison with (\ref{aa1}) and
(\ref{h3s}))
  doubly exponential $L^\infty$-bound (\ref{ff1}) on
 solutions $v(x,t)$. We will then  need in Section \ref{SScal}
  to improve this bound via a
 different scaling approach.

\subsection{Application to a non-divergent equation}

 A similar technique being
applied to the non-divergent diffusion-absorption equation
 \be
 \label{h11}
 v_t= -(-\D)^m v - |v|^{p-1}v \inA,
  \ee
  yields global boundedness of solutions in the subcritical
  Sobolev  range \cite[\S~2]{GW2},
 \be
 \label{h12}
  \tex{
  1 <p <p_{\rm S} = \frac{N+2m}{N-2m} \, .
   }
   \ee
   Further refined applications of the scaling technique  for   (\ref{h11}) are given in
   \cite{ChGal08} (note that Theorem 4.1 proved for $p>p_{\rm S}$ therein in \S~4 applies to
   sufficiently small solutions only).

 \section{Global existence and exponential $L^\infty$-bounds by scaling techniques}
  \label{SScal}

\subsection{Global existence in subcritical range: $C_k$-scaling}

Cf. Theorem \ref{Th.G2}:

\begin{theorem}
 \label{Th.G3}
 The Cauchy problem $(\ref{h1})$, $(\ref{h3})$ has global unique
  classical solution if
 \be
  \label{h10n}
 \tex{
 1<   p  < p_0= 1+ \frac{2(2m-1)}N \quad(m=2l).
  }
   \ee
     \end{theorem}

    Thus, for the standard KSE (\ref{e2.1}) with parameters
     (\ref{KSpar}), in comparison with (\ref{h10n}),
     we manage to skip the first assumption $p \le 3$.



\ssk

\noi{\em Proof.}
  Assume that the local classical solution $v(x,t)$ of the CP
 blows up first time at some finite $t=T$. Then $v(x,T^-)$ is unbounded
 in $L^\infty(\ren)$, otherwise it can be extended as a bounded
 solution on some interval $[T,T+\e]$, with a sufficiently small $\e>0$, by using the
 integral equation (\ref{h5}) with a contractive operator in
 $C(\ren\times[0,\d])$, $\d>0$ small, equipped  with the standard sup-norm.

 Thus, we argue by contradiction and assume that  there exist sequences $\{t_k\}
 \to T^-$, $\{x_k\} \subset \ren$, and $\{C_k\}$ such that
 \be
 \label{seq11}
 \sup_{\ren \times [0,t_k]} \, |v(x,t_k)| = |v(x_k,t_k)| = C_k \to + \infty.
 \ee
Using a modification
  of the rescaling technique in  \cite{GW2},
 we perform the change
 \be
 \label{rvarl}
v_k(x,t) \equiv v(x_k +x, t_k+t) = C_k w_k(y ,s) \whereA x = a_k
y, \quad t = a_k^{2m}s,
 \ee
where $\{a_k\}$ is such that the $L^2$-norm is preserved after
rescaling, i.e.,
 \be
\label{uqq} \|v_k(0)\|_2= \|w_k(0)\|_2 \quad \Longrightarrow \quad
a_k = C_k^{-\frac 2N} \to 0.
 \ee
Therefore, by (\ref{h3s}), for all $s$ for which $w_k(s)$ is
defined,
 \be
 \label{ww1}
  \tex{
 \|w_k(s)\|_2^2
 = \frac 1{a_k^N C_k^2}\, \displaystyle\int v_k^2(x,t) \, {\mathrm d}x
 \le \|v_0\|_2^2 \,\,{\mathrm e}^{\frac T2} \quad
 \mbox{for all} \quad s \in \tex{\big[-\frac{t_k}{a_k^{2m}},
 \frac{T-t_k}{a_k^{2m}}\big)}.
  }
  \ee
 As usual, such a rescaling near blow-up time, in the limit $k \to \iy$ leads to
the so-called  {\em ancient solutions} (i.e., defined for all
$s<0$) in Hamilton's notation \cite{Ham95}, which has been  a
typical technique of
   reaction-diffusion  theory; see various form of its application in
 \cite{SGKM, AMGV}.

 Let, according to (\ref{h10n}),  $p<p_0$.
  Then,  substituting (\ref{rvarl}) into equation (\ref{h1}) yields that
$w_k $ satisfies (as usual, $m=2l$)
 \begin{gather}
 \label{reql}
  (w_k) _s =
  -(-\D)^{2l} w_k + \mu_k (-\D)^l w_k  +  \nu_k \BB_1 |w_k|^{p}  \quad \mbox{in} \quad \ren \times
  \re,\quad \mbox{where} \\
 \mu_k= C_k^{-\frac{2m}N} \to 0, \,\,\, \nu_k = a_k^{2m-1}C_k^{p-1}=
  C_k^{p-1-\frac {2(2m-1)}N} \to 0, \,\, k \to \infty.
 \end{gather}

 We next perform backward shifting in time by
 fixing $s_0>0$ large enough (this is possible in the time-interval in (\ref{ww1}) since $a_k \to 0$ in (\ref{uqq})),
  and setting $\bar w_k (s) = w_k
(s-s_0)$. Then by construction, we have that
 \be
 \label{sc44}
|\bar w_k (s)| \le 1 \andA \|\bar w_k\|_2 \le C \quad \mbox{on}
\,\,\, (0,s_0)
 \ee
 are
uniformly bounded classical solutions of the uniformly parabolic
equation (\ref{Lineq}). By  classic parabolic regularity theory
\cite{EidSys, Fr},
 we have that
 the sequence $\{\bar w_k\}$ is uniformly bounded and  equicontinuous on
 any compact subset from $\ren \times (0,s_0)$. Indeed,
 the necessary uniform gradient bound can be obtained from the
 integral equation for (\ref{reql}), or by other usual regularity methods
 for uniformly parabolic equations. Note that compactness in any
 suitable weaker topologies (see e.g., \cite{Weis0, Weis}) is also
 acceptable, since  passing to the limit we arrive at a weak, and
 hence, classical solution of the limit (simpler) parabolic equation.

 Therefore, by the Ascoli-Arzel\'a theorem,
along a certain subsequence, $\bar w_k (s) \to \bar w (s)$
uniformly on compact subsets from $\ren \times (0,s_0)$. Passing
to the limit in equation (\ref{reql}) and using that both scaling
parameters satisfy $\mu_k, \, \nu_k \to 0$,  yields that $\bar
w(s)$ is a bounded weak solution and, hence, a classical solution
of the Cauchy problem for the linear parabolic equation
(\ref{Lineq}),
 \be
  \label{ww2}
  \bar w_s= -(-\D)^m \bar w, \quad \mbox{with data} \quad
  |\bar w_0| \le 1 \andA
\|\bar w_0\|_2 \le C.
 \ee
  Using the H\"older inequality in the convolution (see the first integral in
  (\ref{h5})) yields
  \be
  \label{ww3}
   \bar w(s_0)= b(s_0)* \bar w_0 \LongA
   |\bar
w(y,s_0)| \le (s_0)^{- \frac N{4m}}\|F\|_2\|\bar w_0\|_2 \ll 1,
   \ee
    for all $s_0 \gg 1$. Hence, the same holds for
$\sup_y|\bar w_k(y,s_0)|$ for $k \gg 1$, from whence comes the
contradiction with the assumption  $\sup_y|v_k(y,s_0)| =1$.

  Thus, $v(x,t)$ does not blow-up and remains bounded for all $t>0$.
 $\qed$

\ssk

The proof is entirely local, so the result holds for the Cauchy
problem as well as for any other homogeneous basic IBVPs for
(\ref{h1}), where the boundary conditions cannot generate blow-up
on the boundary themselves.

\subsection{Exponential $L^\iy$-bound}

We next improve the double exponential $L^\iy$-bound in
(\ref{ff1}):

\begin{theorem}
 \label{Th.G7}
 The global solution of the Cauchy problem $(\ref{h1})$, $(\ref{h3})$
 in the subcritical parameter range $(\ref{h10n})$, $p<p_0$, satisfies
 \be
  \label{h10b}
 \tex{
  \sup_{x \in \ren} \, |v(x,t)| \le C_0 \,{\mathrm e}^{\g_0 t}
  \quad \mbox{for all \,\,$t>0$,\,\, where} \quad \g_0=
  \frac{2m-1}{2N(p_0-p)}\,\,\, \big(>\frac 14 \big).
  }
   \ee
     \end{theorem}

\noi{\em Proof.} We assume that (\ref{seq11}) holds for a monotone
sequence $\{t_k\} \to + \infty$. By (\ref{h3s}), i.e.,
 \be
 \label{hh1}
  \|v(t_k)\|_2 \le \|v_0\|_2 \, {\mathrm e}^{\frac {t_k} 4},
   \ee
  we perform scaling (\ref{rvarl}) by taking into account the
   exponential factor in (\ref{hh1}) targeting a
    uniformly bounded rescaled solution in the sense of (\ref{sc44}).
     This yields (cf.
   (\ref{uqq}))
    \be
    \label{aa21}
     \tex{
    a_k= C_k^{-\frac 2N}\,\, {\mathrm e}^{\frac {t_k}{4N}},
     }
     \ee
   and eventually we arrive at the rescaled equation (\ref{reql}),
   where
  \be
  \label{aa22}
   \tex{
   \mu_k = a_k^{\frac m2}= C_k^{-\frac {m}N}\,\,{\mathrm e}^{\frac {m t_k}{8N}}
    \andA
   \nu_k = C_k^{p-p_0}\,\, {\mathrm e}^{\frac {(2m-1) t_k}{2N}}.
    }
 \ee

Assume now that (\ref{h10b}) does not valid so that, along the
sequence $\{t_k\}$,
 \be
 \label{ddd1}
  C_k = \kappa_k \,{\mathrm e}^{\g_0 t_k}
 \whereA
\kappa_k   \to \infty \asA k \to \infty.
   \ee
One can see that then the rescaled equation (\ref{reql}) contains
the parameters  (recall, $\g_0> \frac 14$)
 \be
 \label{mm1}
  a_k= \kappa_k^{-\frac 2N} \, {\mathrm e}^{- \frac {4\g_0-1}{2N}
  \,t_k} \to 0, \quad \mu_k=a_k^{\frac m2} \to 0, \andA \nu_k =
  \kappa_k^{p-p_0} \to 0.
  \ee
Hence, repeating the arguments of the proof of Theorem
\ref{Th.G3}, we obtain the limit problem (\ref{ww2}), which does
not support the assumed growing behaviour. $\qed$


\subsection{Global existence for the critical exponent $p=p_0$:
$(T-t)$-scaling}

\begin{theorem}
 \label{Th.G6}
 The Cauchy problem $(\ref{h1})$, $(\ref{h3})$ has also global unique
  classical solution in the critical case
 \be
  \label{h10nCr}
 \tex{
  p  = p_0= 1+ \frac{2(2m-1)}N \quad(m=2l).
  }
   \ee
     \end{theorem}

{\bf Remark: on application of the $C_k$-scaling.} For $p=p_0$  in
(\ref{reql}),
 $\mu_k \to 0$, but $\nu_k \equiv 1$, so that passing to the
limit $k \to \infty$, for the limit function $w_k(s) \to w(s)$, we
obtain the KSE in $\ren \times (-\iy,0)$ without the unstable
diffusion-like term:
 \be
 \label{ww2n}
  \tex{
   w_s= -(-\D)^m  w + {\bf B}_1 | w|^p, \quad
  | w(s)| \le 1,\,\,\,
\| w(s)\|_2 \le C, \,\,
 \sup_y|w(y,0)|=1.
 }
 \ee
 Thus, we need to prove that such a solution defined for all $s \le 0$ is nonexistent.

On the one hand, this looks rather reasonable, since
 in the class of uniformly bounded $L^2 \cap H^{2m}$-solutions, the PDE
(\ref{ww2n}) in
 $\ren \times \re_+$ is
a smooth gradient dynamical system  admitting positive definite
Lyapunov function  as in (\ref{h3s}),
 \be
  \label{h3sn}
 \tex{
   \frac 12\,
  \frac{{\mathrm d}}{{\mathrm d}s} \| w(s)\|_2^2=  - \int w(s) (-\D)^m w(s)
  \equiv
  - \|\bar D^m  w(s)\|_2^2  \le 0.
   }
   \ee
   Hence,  the only equilibrium is $0$ that is globally asymptotically
   stable in $L^2(\ren)$, so that, in view of the interior
   parabolic regularity for bounded solutions, for any suitable
   initial data $w_0 \in L^\iy \cap L^2$,
   \be
   \label{jk1}
    w(y,s) \to 0 \asA s \to + \infty \quad \mbox{uniformly}.
    \ee
    On the other hand,
 this is not enough to complete the proof, since
 we  need the convergence
 (\ref{jk1}) to be {\em uniform} with respect to data satisfying
 conditions (\ref{sc44}).
 This will prove the actual nonexistence of a solution to
 (\ref{ww2n}), but is not that straightforward.
 Therefore, for convenience, we choose another, but indeed related, scaling
technique to prove non-blow-up for the critical exponent, which
also emphasize other important aspects of this evolution PDE.

\ssk

\noi{\em Proof of Theorem \ref{Th.G6}.} Thus, assuming first
$L^\iy$-blow-up at $t=T^-$, we perform the time-dependent
$(T-t)$-scaling of the orbit $\{v(\cdot,t), \,t \in (0,T)\}$ of
(\ref{h1}):
 \be
 \label{mm10}
  \tex{
  v(x,t)=(T-t)^{-\a} w(y, \t),\,\,\, \a= \frac {2m-1}{2m(p-1)}, \quad y= \frac x{(T-t)^{1/2m}},
  \quad \t= -\ln(T-t).
   }
   \ee
Without loss of generality, we assume that $x=0$ is a blow-up
point. Then $w$ solves the following exponentially perturbed
equation:
 \be
 \label{mm2}
  \begin{matrix}
 w_\t = {\bf A}(w)+ {\mathrm e}^{-\frac \t 2} (-\D)^l w,
 \quad w_0 \in H^{2m}\cap L^\iy, \ssk\ssk\\
  {\bf A}(w)= -(-\D)^m w- \frac 1{2m} \, y \cdot \n w - \a w+ \BB_1
  |w|^p.
   \end{matrix}
  \ee
Note that the $L^2$-invariance of the scaling (\ref{mm1}):
 \be
 \label{mm3}
  \tex{
  \| v(\cdot, t)\|_2 \equiv \|w(\cdot,\t)\|_2.
   }
   \ee
We next need some auxiliary results.

\begin{proposition}
 \label{Pr.Gl}
 Assume that, along a subsequence $\{\t_k\} \to +\iy$,
  \be
  \label{mm5}
   w(y,\t_k) \to 0 \quad \mbox{uniformly}.
    \ee
 Then $t=T$ is not a blow-up time for $v(x,t)$ (so the singularity  is  removable).
  \end{proposition}

  \noi{\em Proof.} Consider the sequence of solutions
  $\{w_k(y,s)=w(y,\t_k+s)\}$ with vanishing initial data in $L^\iy$
  according to (\ref{mm5}). We now use the good and well-developed
  spectral properties \cite{Eg4} of the linearized operator $\BB^*$ in
  (\ref{mm2}) defined in the weighted space $L^2_{\rho^*}(\ren)$, where
  $\rho^*(y)={\mathrm e}^{-a|y|^\b}$,
  $\b= \frac{2m}{2m-1}$, and  $a>0$ is small enough:
   \be
   \label{mm6}
    \tex{
     \BB^*_\a= -(-\D)^m - \frac 1{2m}\, y \cdot \n- \a I, \quad
   \s(\BB^*_\a)=\big\{ \l_l=-\a- \frac l{2m}, \,\,\,
   l=0,1,2,...\big\},
   }
   \ee
  with   eigenfunctions (generalized
   Hermite polynomials) that compose a complete and closed set.
   Therefore, according to classic asymptotic parabolic theory
   (see e.g., \cite{Lun}), we conclude that for any sufficiently
   large $k$,
    \be
    \label{mm7}
     \tex{
     w_k(y,s) \sim O\big({\mathrm e}^{- \a s}\big), \,\,\, s \gg 1
     \LongA  w(y,\t) \sim O\big({\mathrm e}^{- \a \t}\big), \,\,\, \t \gg
     1.
     }
     \ee
Overall, taking into account scaling (\ref{mm10}), this yields:
 \be
 \label{mm8}
 \tex{
 v(x,t) \sim (T-t)^{-\a} O\big({\mathrm e}^{- \a \t}\big)= O(1) \asA
 t \to T^-,
 }
 \ee
 so that $v(x,t)$ is uniformly bounded at $t=T$. $\qed$

 \ssk

\begin{proposition}
 \label{Pr.Gl2}
 There must exist a subsequence $\{\t_k\} \to +\iy$ such that
  \be
  \label{mm5N}
C_k \equiv   \|w(\cdot,\t_k)\|_\iy \to \iy \quad \mbox{as \,\, $k
\to +\iy$}.
    \ee
  \end{proposition}


  \noi{\em Proof.} Assume for contradiction that
  \be
  \label{mm9}
  w(y,\t) \quad \mbox{is uniformly bounded in $\ren \times
  \re_+$}.
   \ee
Then, similar to (\ref{h3sn}), we get for uniformly bounded and
smooth solution $w(y,\t)$ that
\be
  \label{mm11}
 \tex{
   \frac 12\,
  \frac{{\mathrm d}}{{\mathrm d}\t} \| w(\t)\|_2^2=
  - \|\bar D^m  w(\t)\|_2^2 +{\mathrm e}^{-\frac \t 2}\|\bar D^l
  w(\t)\|_2^2.
   }
   \ee
 Therefore, in order to avoid the $L^2$-vanishing:
  \be
  \label{mm12}
\| w(\t)\|_2^2 \to 0 \asA \t \to +\iy,
 \ee
 which by the interior parabolic regularity would imply
 (\ref{mm5}) as $\t \to +\iy$ uniformly and hence no blow-up, one
 has to have that
  \be
  \label{mm13}
  \tex{
   \int^\iy \|\bar D^m  w(\t)\|_2^2 \,\,{\mathrm d}\t < \iy, \quad
   \mbox{i.e.,} \quad
\|\bar D^m  w(\t)\|_2^2 \in L^1((1,\iy)).
 }
 \ee
Obviously, the convergence of the integral in (\ref{mm13}) implies
that
 \be
 \label{mm14}
 \tex{
 \bar
D^m w(\cdot, \t + s) \rightharpoonup 0 \,\,\, \mbox{as $\t \to
\iy$ weakly in $L^2_{\rm loc}(\re_+;L^2)$}. }
 \ee
Indeed, one can see that, by the H\"older inequality, for any
$\chi \in C_0^\iy(\ren \times \re_+)$,
 \be
 \label{mm15}
   \begin{matrix}
  \big|-\iint (-\D)^m w(\t+s) \chi(s)\big| \equiv  \big|\iint (\bar D^m)^2 w(\t+s) \chi(s)\big|
  \qquad
  \ssk\ssk\\
= \big|\iint \bar D^m w(\t+s) \bar D^m \chi(s)\big| \le \big(
\iint( \bar D^m w(\t+s))^2 \big)^{\frac 12} \big( \iint (\bar D^m
\chi(s))^2 \big)^{\frac 12} \to 0 \qquad
\end{matrix}
 \ee
 as $\t \to \iy$. Fixing a sequence $\{\t_k\} \to \iy$ and passing to the
limit as $\t_k + s \to \iy$, we conclude that, in view of the
remaining interior parabolic regularity (recall that, regardless
the ``artificial degeneracy" (\ref{mm14}), equation (\ref{mm2}) is
uniformly parabolic), some equicontinuous subsequence
$\{w(\t_k+s)\} \to \hat w(s)$ uniformly,
 where $\hat w(y,s)$ is a smooth solution
of  the first-order Hamilton--Jacobi equation
 \be
 \label{m16}
  \tex{
 \hat w_s= - \frac 1{2m} \, y \cdot \n \hat w- \a \hat w+
 \BB_1|\hat w|^p, \quad \hat w_0 \in L^2 \cap L^\iy.
 }
  \ee
However, by classic theory of conservation laws \cite{Tay} (and,
actually, by the standard comparison), all such solutions of
(\ref{m16}) are exponentially decaying:
 \be
 \label{m17}
 \hat w(y,s) = O\big({\mathrm e}^{-\a s}\big) \asA s \to \iy.
 \ee
 This implies that there exists always a moment $\t_0=\t_k+s_0$, with large enough $k$, $\t_k$,
 and $s_0$ such
 that $ w(y,\t_0)$ is arbitrarily
uniformly. As in the proof of Proposition \ref{Pr.Gl}, overall,
this implies no blow-up for $v(x,t)$ at $t=T$. $\qed$

\ssk

Finally, the result on global existence is completed by the
following:

\begin{proposition}
 \label{Pr.Gl3}
 The Type II\footnote{The terms ``Type I, II" were borrowed from Hamilton
 \cite{Ham95}, where Type II is
  also called {\em slow} blow-up. In reaction-diffusion theory \cite{SGKM, AMGV},
  Type I blow-up is usually called of self-similar rate, while Type II is referred to as fast and
   non-self-similar.} blow-up solutions satisfying $(\ref{mm5N})$ do not
 exist.
  \end{proposition}

  \noi{\em Proof.} We now apply in (\ref{mm2}) the $C_k$-scaling with $\{C_k\}$
  given in (\ref{mm5N}), i.e., as in (\ref{rvarl}).
  For $p=p_0$,
  this gives the following equation
 this gives the following equation for the rescaled sequence:
 \be
 \label{sg1}
  w_k(y,\t) \equiv w(y_k+a_k z, \t_k+a_k^{2m}s)= C_k \hat w_k(z,s),
  \,\,\,\, \mbox{where \,$\hat w_k(z,s)$\, solve}
  \ee
  \be
 \label{mm2NN}
 \tex{
\hat w_s = -(-\D)^m \hat  w + \BB_1
  |\hat w|^p - a_k^{2m}\big( \frac 1{2m} \, z \cdot \n \hat w+ \a
  \hat w\big) + a_k^m
  {\mathrm e}^{-\frac {\t_k+a_k^{2m}s} 2} (-\D)^l \hat w.
  }
  \ee
   Here, at $s=0$,  $\hat w_0 \in L^2\cap L^\iy$ and $a_k=C_k^{- 2/N}$.
This is a uniformly parabolic equation with asymptotically small
perturbations, so that on passage to the limit $k \to \iy$ for
this set $\{\hat w_k\}$ of  smooth solutions, one has to have in
the limit that, along a subsequence,
 \be
 \label{mm29}
   \begin{matrix}
w_k=w(y_k+a_k z,\t_k+a_k^{2m} s) \to W(z,s) \whereA \quad
 \ssk\ssk\ssk\\
 W_s=
 -(-\D)^m W + \BB_1
  | W|^p\inB \ren\times \re, \,\, W \in L^2 \cap L^\iy, \,\,
  \|W(z,0)\|_\iy=1.
  \quad
 \end{matrix}
   \ee
 Note that $W(z,s) \not \equiv 0$ is an ancient solution, which is defined for
 all $s \le 0$. At the same time, by construction, it is also a
 {\em future solution}, which must defined for all $s>0$. Indeed,
 one can see that if $W(s)$ blows up at some finite $s=S^->0$,
 this would contradict the Type II solution $w(y,\t)$ is globally
 defined for all $\t > 0$.
  Thus, by scaling of the Type II blow-up orbit (\ref{mm5N}), we arrive at the problem (\ref{mm29}),
  which defines:
 \be
 \label{nn1}
 \mbox{$\{$\underline{heteroclinic} solution $W(z,s) \not =0\}$ = $\{$\underline{ancient} for
 $s<0\}\,\,\cup\,\,\{$\underline{future} for $s>0\}$}.
 \ee

Let us more clearly specify the necessary properties of this
mysterious and hypothetical  heteroclinic solution $W$.
   Bearing in mind the possibility of multiplying (\ref{mm29}) by
   $W$
and integrating by parts, which are again guaranteed by the
convergence demand as in  (\ref{mm13}), we have the standard
identity:
 \be
 \label{mm66}
 \tex{
   \frac 12\,
  \frac{{\mathrm d}}{{\mathrm d}s} \| W(s)\|_2^2=
  - \|\bar D^m  W(s)\|_2^2 \le 0,\,\, s \in \re \LongA
   \int\limits_{-\iy}^{+\iy}  \|\bar D^m  W(s)\|_2^2\, {\mathrm d}s <
   \iy.
   }
   \ee
   Hence, in the given class $L^2 \cap L^\iy$ classical solutions, (\ref{mm29})
   is a sufficiently smooth gradient system, and the only
   equilibrium that can be approached by such bounded orbits is
   $W=0$.
    By passing to the
   limit $s \to \pm \iy$, we then obtain that uniformly
    \be
    \label{nn2}
    W(s) \to 0 \asA s \to \pm \iy.
    \ee
This actually means that
 \be
 \label{nn3}
 \mbox{the heteroclinic solutions $W(z,s) \not \equiv 0$ is a
 \underline{homoclinic of 0 orbit}}.
 \ee
We next easily prove that

\begin{proposition}
 \label{Pr.Anc}
 A nontrivial solution $W \not =0$ of $(\ref{mm29})$ does not
 exist.
  \end{proposition}

  \noi{\em Proof.} It suffices to use the ``ancient" part of the
  definition (\ref{nn1}). In view of the integral convergence in
  (\ref{mm66}) at $s=-\iy$, similar to (\ref{mm5}) we have that
  along any monotone sequence ${s_k} \to - \iy$, we have that, for
  $k \gg 1$,
  $W(s_k+s) \approx V(s)$ uniformly on compact subsets in $\ren \times
  \re$, where $V$ solves the  Hamilton-Jacobi equation (a
  conservation law)
   \be
   \label{nn5}
   V_s=\BB_1|V|^p \forA s>0, \quad V_0=W(s_k) \in L^2\cap L^\iy.
    \ee
Since in the class of smooth solutions, $V(s)$ decays fast (see
e.g., \cite{Sm, Tay}), we have that
 \be
 \label{nn7}
 \|V(0)\|_\iy \gg \|V(s)\|_\iy \forA s \gg 1,
  \ee
  so that the same is true for $W(s_k+s)$ provided that $k \gg 1$.
In fact, (\ref{nn7}) means that using proper theory of
conservation laws (\ref{nn5}) implies that a nontrivial $L^2\cap
L^\iy$ ancient solution does not exists for such gradient system,
since by the convergence in (\ref{mm66}), equation (\ref{mm29})
does not have a sufficient mechanism for growth of solutions from
$\|W(-\iy)\|_\iy=0$ to $\|W(0)\|_\iy=1$. To justify this more
clearly, given an ancient solution $W(z,s)$ of (\ref{mm29}), we
perform the standard scaling:
 \be
 \label{WW1}
  \tex{
  W_\l(z,s)=\l^{-\a} W\big( \frac z{\l^{1/2m}}, \frac s \l\big),
  \quad \l >0,
  }
  \ee
  where $W_\l$ is also an ancient solution of (\ref{mm29}) for any $\l >0$.
  We now pass to the limit $\l \to 0^+$ by using these simple
  properties:
   \be
   \label{WW2}
    \begin{matrix}
 \|W_\l(0)\|_\iy = \l^{-\a} \to +\iy, \quad
    \|W_\l(s)\|_2^2 = \|W(\frac s \l)\|_2 \to c_0>0,
     \ssk\ssk\ssk \\
     \andA \frac 12\, \frac{\mathrm d}{{\mathrm d}s}\,
     \|W_\l(s)\|_2^2 =-\l \| \bar D^m W(\frac s \l)\|_2^2 \to 0.
  \end{matrix}
      \ee
Therefore, in the limit $\l \to 0$, we have to have  a nontrivial
$L^2$-solution $V(z,s)$ of (\ref{nn5}), which blows up as $s \to
0^-$. Obviously, this contradicts the Maximum Principle for this
first-order conservation law, \cite[Ch.~16]{Tay}.

 In other words, we have that, for the dynamical system in (\ref{mm29}),
 \be
 \label{cc1}
\mbox{the unstable manifold of the origin 0 is empty}.
  \quad \qed
  \ee



This also completes the proof of Theorem \ref{Th.G6}. $\qed$

\subsection{On uniform bounds in other KS-type models}


 For the non-divergent equation (\ref{h11}), a similar
scaling techniques yields that \cite{GW2}, in the Sobolev range
(\ref{h12}), solutions are uniformly bounded, i.e.,
 \be
  \label{bb1}
 |v(x,t)| \le C \inA.
 \ee


For the divergent mKS-type equations {\em without} the unstable
backward diffusion term,
 \be
 \label{h1N}
 v_t=-(-\D)^m v +  \BB_1 (|v|^p) \inA  \quad (m \ge 2),
  \ee
 we easily extend the result as follows:

 \begin{theorem}
 \label{Th.G4}
 The Cauchy problem $(\ref{h1N})$, $(\ref{h3})$ in the parameter range $(\ref{h10n})$
  has a global unique
  classical solution,
 which is uniformly bounded, i.e., $(\ref{bb1})$ holds.
     \end{theorem}

\noi{\em Proof.} Again, it suffices to consider the case when
 (\ref{seq11}) holds for some sequence $\{t_k\} \to +\infty$. Then
 the same proof leads to the contradiction with the hypothesis
 that $\{C_k \} \to + \infty$, and hence $v(x,t)$ is global and
 uniformly bounded. $\qed$

\subsection{On a generalization: higher-order nonlinear dispersion}

In order to extend application of our final approach,
 we next briefly discuss  mKS-type equations
with a {\em third-order} (also {\em odd}) nonlinear perturbation
(dispersion) of the form
\be
 \label{h1N3}
 v_t=-(-\D)^m v - \D\BB_1 (|v|^p) \inA  \quad (m \ge 2).
  \ee
Writing this PDE in a pseudo-parabolic form,
 \be
 \label{h1N33}
 P v_t=(-\D)^{m-1} v +  \BB_1 (|v|^p) \whereA P=(-\D)^{-1}>0,
  \ee
 and multiplying by $v$ in $L^2(\ren)$,
  we observe that, instead of a uniform $L^2$-bound, we are given
  an {\em a priori} $H^{-1}$-bound: for uniformly bounded data
   \be
   \label{vv1}
   v_0 \in L^\infty(\ren) \cap
  H^{2m}(\ren),
  \ee
   the following  holds:
   \be
   \label{bb21}
   \|v(t)\|_{-1} \le \|v_0\|_{-1} \forA t>0.
    \ee
    Here, for simplicity, we assume that $v_0(x)$ has also
    exponential decay at infinity, so $v(x,t)$ does for $t>0$.
In (\ref{h1N33}) and later on, we define $(-\D)^{-1} w=g$ in a
standard manner:
 \be
 \label{s1}
 \D g = - w \inB \ren, \quad g(x) \to 0 \asA x \to \infty.
  \ee
For the  solvability of this problem we shall  always assumes that

  \be
   \label{s2}
  \tex{\displaystyle
   \int\limits_{\ren} w\,{\mathrm d}x=0.
   }
   \ee
Clearly this property holds for the divergent equation
 (\ref{h1N3}) with exponentially decaying solutions.

 \begin{theorem}
 \label{Th.G5}
 A unique solution $v(\cdot,t) \in H^{-1}(\ren)$ for $t>0$ of
 the Cauchy problem $(\ref{h1N3})$, $(\ref{h3})$, $(\ref{vv1})$
  is uniformly bounded, i.e., $(\ref{bb1})$ holds, in the
  parameter range
   \be
   \label{ra1}
 \tex{
 1<p < p_0=1 +\frac{2(2m-3)}{N+2}.
 }
   \ee
     \end{theorem}

\noi{\em Proof.} Local existence of classic solutions for
(\ref{h1N3}) also follows from the equivalent integral equation
such as (\ref{1.2N}), with no second unstable term by replacing
 $$
 \BB_1^* \mapsto -\BB_1^* \D.
  $$
   One can see that we still obtain a
locally integrable in $t>0$ kernel that defines the operator being
a contraction in $C(\ren \times [0,\d])$, with sup-norm.

To prove global and uniform boundedness,
 we use the
same scaling as in the proof of Theorem \ref{Th.G3}, where,
instead of (\ref{uqq}), keeping the $H^{-1}$-norm of $w_k$ yields
 \be
\label{uqq3}
 \|v_k(t_k)\|_{-1}= \|w_k(0)\|_{-1} \quad \Longrightarrow
\quad a_k = C_k^{-\frac 2{N+2}}.
 \ee
Then, for $p<p_0$, we arrive at the corresponding equation such as
(\ref{reql}), where
 \be
  \label{jj1}
   \tex{
   \nu_k=a_k^{2m-3} C_k^{p-1}= C_k ^{
   p-1-\frac{2(2m-3)}{N+2}} \to 0 \quad (\mu_k=0)
   }
   \ee
as $k \to \infty$. By passing to the limit $k \to \infty$, we
again arrive at the purely poly-harmonic flow  as in (\ref{ww2})
that cannot support the necessary properties of the sequence. As
usual, this analysis applies twice: (i) $\{t_k\} \to T^-<\infty$,
to prove that no blow-up occurs, and (ii) $\{t_k\} \to +\infty$,
to prove that the solution is uniformly bounded. $\qed$

\ssk

The critical case $p=p_0$ can be also studied in similar lines as
above. Note in addition that  0 is the unique globally
asymptotically stable equilibrium of the smooth gradient system
(\ref{h1N3}) in the corresponding class of regular solutions that
are uniformly bounded in $L^\iy \cap H^{-1}$.

\ssk

Further generalizations  including  odd higher-order nonlinear
dispersion operators are straightforward.

\subsection{On blow-up in divergent models}

Phenomena of finite-time  blow-up in such semilinear models with
divergent operators is most well-known for {\em unstable limit
Cahn--Hilliard equation}
\be
\label{GPP}
    u_t = -\Delta^2 u -\D(|u|^{p-1}u) \quad \mbox{in} \quad \ren \times
    \re_+\quad (p>1).
\ee Blow-up solutions have the standard self-similar form
\be
\label{ResVars}
 u(x,t) = (T-t)^{-\frac 1{2(p-1)}}  f(y),
\quad y = {x}/(T-t)^{\frac 14},
 \ee
where $f$ solves the elliptic equation
\be
\label{RescOp}
  - \Delta^2 f -\D(|f|^{p-1}f) -
\mbox{$\frac{1}{4}$}\, y\cdot \nabla f - \mbox{$
\frac{1}{2(p-1)}$}\, f=0 \inB \ren.
 \ee
This equation admits a complicated set of solutions. For instance,
for $N=1$ and $p=3$, it has a countable set of different profiles
$f$ that describe various types blow-up patterns (\ref{ResVars}).
See \cite{EGW1} for earlier related references and  further
results. The case $p=2$ was earlier considered in \cite{BB95}.
Concerning {\em Leray's scenario} of similarity blow-up in
(\ref{GPP})  in both limits $t \to T^\mp$, see \cite{GalJMP} and
references therein.

\ssk

We are not aware of  reliable traces of standard blow-up for
KS-type PDEs such as (\ref{h1}) or (\ref{h1N3}) with odd-order
nonlinear dispersion terms in the supercritical range $p>p_0$.
Therefore, we do not know whether conditions (\ref{h10n}) and
(\ref{ra1}) of global solvability reflect the actual evolution
properties of the PDEs under consideration, or are sometimes
purely technical. In other words, one then faces  another
fundamental problem on construction of blow-up patterns for
$L^2$-bounded solutions.  Regardless a  sufficiently strong
progress in understanding of formation of blow-up singularities in
various nonlinear PDEs achieved in the last twenty five years (see
references and results in the monographs \cite{GalGeom, AMGV,
MitPoh, SGKM, QSupl}), for higher-order equations, there are still
several fundamental open problems in identifying admitted
structures of blow-up patterns.


 \section{On $L^\infty$-bounds for the
 Navier--Stokes equations in $\ren$ and well-posed Burnett equations}
 \label{SNS}

\subsection{ A classical fluid model in $\ren$}

Consider the Navier--Stokes equations (\ref{NS1}) with given
bounded $L^p$-data $\vv_0$. In order to apply our scaling
argument, we use the fact that  a classical bounded solution
$\vv(x,t)$ can be locally extended by using the integral equation
that is similar to \eqref{h5}.
 Existence of such a local semigroup of smooth bounded solutions for
 the NSEs is well known; see Majda--Bertozzi \cite{Maj02}.

Let us present some comments that will be useful for the Burnett
equations \eqref{NS1m} with $m \ge 2$.
 Taking the fundamental solution \eqref{1.3R} for $m=1$ with
  the rescaled Gaussian
  \be
  \label{G1}
   \tex{
 F(y)=( 4 \pi)^{-\frac N2}{\mathrm e}^{- \frac{|y|^2}4},
  }
  \ee
 we consider (\ref{NS1}) with $\hh=-\n p$ as a system
  \be
  \label{G2}
  \left\{
   \begin{matrix}
  \vv(t)= b(t) * \vv_0- \displaystyle\int\limits_0^t b(t-s)*[(\vv\cdot \n)\vv](s) \,
  {\mathrm d}s +\displaystyle\int\limits_0^t b(t-s)*\hh(s) \,
  {\mathrm d}s,\qquad\quad\,\, \\
  \di b(t) * \vv_0-\displaystyle\int\limits_0^t \di b(t-s)*[(\vv\cdot \n)\vv](s) \,
  {\mathrm d}s +\displaystyle\int\limits_0^t \di b(t-s)*\hh(s) \,
  {\mathrm d}s=0.
   \end{matrix}
   \right.
\ee
 As usual, the second equation in (\ref{G2}) is  the one for
the pressure corresponding to the  solenoidal vector field $\vv$.
 Observe that, due to the exponential decay of the Gaussian (\ref{G1}), the first equation
  contains the operator in $\vv$ that
 contractive in a bounded closed subset $M_\d$ of $C([0,\d], C^1(\ren))$,
 where $\d>0$
 sufficiently small, with the sup-norm.
 Indeed, assuming that in $M_\d$,
  $$
  | \vv| \le C \andA | D \vv | \le C,
  $$
  we can use the possibility of differentiating in $x$  the
  equation once to control $D \vv$. For fixed vectors $\vv_0$ and
  $\hh(s)$, the contractivity of the principal ``convective" operator
   \be
   \label{op11}
    \tex{
    {\mathcal N}(\vv)=
\displaystyle\int\limits_0^t b(t-s)*[(\vv\cdot \n)\vv](s) \,
  {\mathrm d}s
  }
   \ee
    in $M_\d$ for small $\d>0$ is then
   straightforward.
   Note that
   the standard Picard iteration scheme for (\ref{G2}) can be put
   into the probability framework by using the linear
   semigroup associated with 3D Brownian motion; see a survey \cite{Way05}
   for these and other details related to the Navier--Stokes equations.
 We stop at this moment discussing the local interior regularity
 theory for the Navier--Stokes equations and refer to
 \cite{Can05, Gust07, Maj02} as a guide for detailed developments in this direction.

As a standard classic alternative way for local applications,  the
pressure is excluded from the NSEs
 \be
  \label{ww1N}
     \vv_t={\bf H}( \vv) \equiv  \D \vv
 -  \mathbb{P}\,(\vv
\cdot \n) \vv, \,\,\,\,
   \mbox{where} \quad \mathbb{P} \uu=\uu - \n \D^{-1}(\n \cdot
   \uu)
  \ee
  is the
  Leray--Hopf projector onto the solenoidal vector
    field.
 Using the fundamental
 solution of $\D$ in $\ren$, $N \ge 3$ ($\s_N$ is the surface area of the unit ball $B_1 \subset \ren$)
  \be
  \label{FF55}
   \tex{
   b_N(y)= - \frac 1{(N-2)\s_N} \, \frac 1{|y|^{N-2}} \whereA \s_N= \frac { 2
   \pi^{ N/2}}{\Gamma(\frac N2)},
    }
    \ee
the operator in \eqref{ww1N} is written in the form of Leray's
formulation
\cite[p.~32]{Maj02} 
 \be
 \label{HHH21}
  \begin{matrix}
{\bf H}( \vv) \equiv  \D \vv - (\vv \cdot \n ) \vv + C_3
\displaystyle\int\limits_{\re^3} \tex{\frac {y-z}{|y-z|^3}}\,\,
{\rm tr} (\n  \vv (z,\t))^2\, {\mathrm d}z, \ssk \\
 \mbox{where} \quad
 {\rm tr} (\n
\vv (z,\t))^2= \sum_{(i,j)} \,   v_{z_j}^i  v_{z_i}^j \andA C_N=
\frac 1{\s_N}>0.
 \end{matrix}
 \ee
 Equivalently,
 the nonlocal parabolic equation \eqref{ww1N} is  known to induce a
 local semigroup of smooth solutions.

As usual, our intention is to show
 that local sufficiently smooth solutions cannot blow-up
 under a certain
$L^p$-type constraint. We next clarify the conditions that prevent
finite-time blow-up of solutions in $L^\infty$ and exclude Leray'
similarity \eqref{NS2} or any other.

It is a classic matter  that the $L^2$-norm is natural for
(\ref{NS1}), since after multiplication by $\vv$, the convective
and pressure terms vanish on sufficiently smooth functions
$\vv(x,t)$ with fast decay at infinity,
 \be
 \label{na1}
  \langle (\vv \cdot \n)\vv, \vv \rangle =0
  \andA -\langle \n p, \vv \rangle = \langle p, \n \cdot \vv
  \rangle=0.
   \ee
Therefore, on  smooth solutions,
 \be
 \label{ss1}
  \tex{
   \frac 12\, \frac{\mathrm d}{{\mathrm d}t} \, \|\vv(t)\|_2^2 = -
   \| D \vv(t)\|_2^2 \LongA \|\vv(t)\|_2 + 2 \displaystyle\int\limits_0^t \| D \vv(t)\|_2^2
    \, {\mathrm d}t \le \|\vv_0\|_2, \,\,\, t
   \ge 0. }
    \ee
Actually, the estimate in (\ref{ss1}) is the energy inequality
 for  {\em Leray--Hopf weak solutions} of (\ref{NS1}); see e.g.,
 \cite{Qio07} and references therein.

Nevertheless,
an $L^2$-bound is not sufficient to control non-blowing up
property of solutions, and this is   the origin of extensive
mathematical research in the last fifty years. Recall that
 global weak solutions of (\ref{NS1}) satisfying
  $$
  u \in L^\iy(\re_+;L^2(\re^3)) \cap L^2(\re_+;H^1(\re^3))
  $$
  were constructed by Leray \cite{Ler34}, and Hopf \cite{Hopf51}
  in 1951.

\subsection{Application of blow-up $C_k$-scaling in the supercritical case $p>N$}

Namely, as a first simple constraint, which is  inherited from our
previous study of the KS-type equations, we assume that
 \be
 \label{rr1}
 \| \vv(t)\|_p \le C \forA t \ge 0 \quad( p > 2).
  \ee
Then, seeking $L^\iy$-bound and hence assuming \eqref{seq11}, we
perform the scaling (\ref{rvarl}), where we impose the
preservation of the $L^p$-norm of the rescaled function, i.e.,
 \be
  \label{rr2}
  \|\vv_k\|_p=\|\ww_k\|_p \LongA  a_k= C_k^{-\frac pN} \to 0,
  \andA
  b_k=a_k^2.
   \ee

As usual, we next perform passage to the limit in the NSEs.
This can be done in the framework of the original model
\eqref{NS1}, as well as of the nonlocal parabolic representation
\eqref{ww1N}, which we actually do. Note that passage to the limit
in the integral term causes no difficulty for our sequences of
uniformly bounded smooth solutions
 $\{\ww_k\}$, where  Ascoli-Arzel\'a classic theorem
 \cite[Ch.~2]{KolF} applies.

\ssk

We then obtain the following rescaled equations for
$\ww=\ww_k(y,s)$:
  \be
  \label{NS1n}
 \ww_s + \nu_k \, \mathbb{P} (\ww \cdot \n)\ww=
 \D \ww, \quad
   \mbox{where} \quad \nu_k=C_k^{1- \frac pN} \to 0 \forA p>N.
    \ee
  Next, after time shifting, $s \mapsto s-s_0$,
 the solutions and data satisfy (cf. \eqref{sc44})
  \be
  \label{rr4}
  | \bar w_k| \le 1 \andA \|\bar w_k\|_p \le C.
   \ee
 By  regularity for uniformly bounded sequence
 $\{\bar \ww_k(s)\}$, we conclude that there exists its partial limit
$\bar \ww(s)$ satisfying the {\em solenoidal heat equation}
 \be
 \label{NS1n1Heat}
\bar \ww_s  = \D \bar\ww.
 \ee
Of course, this is an equivalent form of writing
 the nonstationary linear
``convectionless" Stokes system,
  \be
  \label{NS1n1}
 \bar \ww_s  =- \n \bar q + \D \bar\ww, \quad \di \bar\ww=0
  \inA,
  \ee
 with data as in (\ref{rr4}).
 It is an exercise to check, by using the integral equation
 with the kernel as in (\ref{h5}), that this problem for (\ref{NS1n1Heat}) or (\ref{NS1n1})
 admits a solution that exhibits a power decay for $s \gg 1$. Namely,
 an estimate similar to (\ref{ww3}) can be obtained with a
 different exponent of the decay rate.

  In other words, under the assumption (\ref{rr1}), the
 Navier--Stokes system does not have a mechanism to create any
 $L^\infty$ blow-up singularities. We thus arrive at the
 following:

  \begin{proposition}
 \label{Th.Ex}
 Under the given assumptions,  $L^p$-solutions, with $p > N>2$ in $\eqref{rr1}$
of the Navier--Stokes equations $(\ref{NS1})$ do not blow-up, and,
moreover,
 are uniformly
bounded:
 \be
 \label{ka1}
 |\vv(x,t)| \le C \inA.
  \ee
  \end{proposition}

The last bound (\ref{ka1}) is proved as in Theorem \ref{Th.G7} by
assuming that $\{t_k\} \to + \infty$.

\ssk

One can change the functional space in estimate (\ref{rr1}) to get
$L^\iy$-bound.


\subsection{Critical case $p=N=2$: an example of application of the $(T-t)$-scaling}
 \label{S6.3}

In the critical case $p=N$, we have that $\nu_k=1$ in
(\ref{NS1n}), so that in the limit $k \to +\infty$ we arrive at
the same Navier--Stokes equations, but now with {\em uniformly
bounded data and solutions} in both $L^3(\ren)$ and
$L^\infty(\ren)$.
 Moreover, by passing to the limit $s_0=s_{0j} \to +\iy$,  we actually deal with the
  following class of solutions:
  \be
  \label{abs1}
   \mbox{$\bar \ww(s) \in L^\iy\cap L^N$ for all $s \ge 0$: $|\ww(s)| \le 1$, \,\,
   $\|\ww(s)\|_N \le C$.}
    \ee
 In other words, using the scaling in the critical case $p=N$, we
 eventually get into the special class of solutions (\ref{abs1}),
 so a key restriction
 of the Navier--Stokes equations is achieved.
Obviously, no blow-up and other singularities are available in the
class (\ref{abs1}).
Thus,  in the critical case,  in the class (\ref{abs1}), the
Navier--Stokes equations induce a smooth gradient dynamical system
with the positive Lyapunov function as in (\ref{ss1})
 which is strictly monotone on such nontrivial solutions.
 Hence, this  admits
 the
unique globally asymptotically stable trivial equilibrium $0$.

However, for $p=N$, proving non blow-up of solutions leads to a
hard problem of nonexistence of suitable ancient solutions.
Similar to the previous KS problem, we demonstrate an example of
application of the $(T-t)$-scaling to get the result in the
critical case $p=N=2$ of the obvious particular interest.

We begin with Leray's blow-up scaling \cite{Ler34} for
(\ref{ww1N}) by setting
 \be
 \label{w1}
 \tex{
 \vv(x,t)= \frac 1{\sqrt{T-t}}\,  \ww(y,\t), \quad y= \frac
 x{\sqrt{T-t}}, \quad \t=-\ln(T-t),
  }
  \ee
  to get the following rescaled equation:
   \be
   \label{w2}
    \tex{
     \ww_\t+ {\mathbb P}( \ww \cdot \n)  \ww= \BB^*
    \ww \whereA \BB^*= \D - \frac 12\, y \cdot \n - \frac 12\,I.
     }
     \ee
Here, $\BB^* $ is the adjoint Hermite operator with the discrete
spectrum
 \be
 \label{w3}
  \tex{
  \s(\BB^*)= \big\{ \l_k=-\frac 12- \frac k2, \,\,k=0,1,2,...\big\}
  }
  \ee
  and a complete and closed set of eigenfunctions being finite
  solenoidal Hermite polynomials; see details in \cite[Append.~A]{Gal02A} and
  \cite[\S~2]{GalNSE}, where all the aspects of the functional
  settings of $\BB^*$ in the weighted space $L^2_{\rho^*}$, $\rho^*(y)={\mathrm
  e}^{-|y|^2/4}$, can be found. Note that scaling (\ref{w1})
  implies the following behaviour of the rescaled solution:
   \be
   \label{ww87}
    \|\ww(\t)\|_2^2 \sim {\mathrm e}^{\frac{N-2}2\, \t}\big|_{N=2}
    =O(1) \asA \t \to +\iy.
      \ee


We next follow the same scheme as in Section \ref{SScal}:

\noi{\bf (i)} First, the $L^2$-conservation such as (\ref{mm3})
holds for $\ww(y,\t)$; see (\ref{ww87}).

\noi{\bf (ii)} Second, assuming (\ref{mm5}) yields by using the
spectral properties in (\ref{w3}) and scaling (\ref{w1}) (cf.
(\ref{mm8})) the non-blow-up in $L^\iy$ at $t=T$ of $\vv(x,t)$.

\noi{\bf (iii)} Thirdly, we arrive at (\ref{mm5N}) for
$\ww(y,\t)$, since its uniform  boundedness will lead to
$L^\iy$-decay by the identity
 \be
 \label{w6}
  \tex{
  \frac 12 \,\frac{\mathrm d}{{\mathrm d}\t} \, \|\ww(\t)\|_2^2 =
  -
  \int |\n \ww|^2 <0.
   }
   \ee
   This $\n$-property of the dynamical system means that, in the
   class of bounded smooth solutions, the only equilibrium is
   trivial $\ww=0$, so that the uniform stabilization to it
   again guarantees non-blow-up.

\noi{\bf (iv)} Fourth, assuming (\ref{mm5N}), we $C_k$-scale such
an orbit to get a convergence to a  nontrivial 0-homoclinic
solution as in (\ref{mm29}):
 \be
 \label{w7}
  \tex{
  \WW_s+
  {\mathbb P}(\WW \cdot \n)\WW= \D \WW \inB \ren \times \re,
   \quad \int\limits_{-\iy}^\iy \|\n \WW(s)\|_2^2\, {\mathrm d}s <
   \iy,
   }
   \ee
   where, thus, $\WW(s) \in L^2 \cap L^\iy$, and, by (\ref{w6}),
    $\WW(s) \in H^1$ for all $s \in \re$. Otherwise, if $\WW \not \in H^1$, then
    (\ref{w6}) will guarantee fast decay of the $L^2$-norm of
    $\ww(\cdot,\t)$ and hence $\vv(\cdot, t)$ as $t \to T^-$, which is
    contradictory. Finally, using the converging integral in
    (\ref{w7}) and passing to the limit as $s_k+s \to -\iy$
    leads to smooth small solutions of the Euler equations (EEs):
     \be
     \label{VV1}
     \VV_s+ {\mathbb P}(\VV \cdot \n)\VV=0, \quad \VV(s) \in L^2
     \cap L^\iy.
     \ee
 Then the void conclusion similar to (\ref{cc1}) remains valid,
 provided that smooth $L^2$-solutions of the EEs (\ref{VV1}) decay in
 $L^\iy$ sufficiently fast, which is true since this is a smooth
 gradient system (recall that we are obliged to treat the simpler critical case $p=N=2$ only);
  see \cite[Ch.~17]{Tay} and surveys \cite{Bar07, Con07} for further details.
This shows a potential correct way to treat the relations between
singularities in the  NSEs and EEs, where the absence of some of
those for the latter ones implies the same for the former.


Of course, our analysis of global existence of classical bounded
solutions in the critical case $p=N=2$  just reflects the classic
Leray's (1933, the CP)
 \cite{Ler33N2}
 and Ladyzhenskaya (1958, IBVPs)
 \cite{Lad58, Lad61} (see also \cite{Lad70})
existence-uniqueness results for $N=2$. In this connection, it is
worth mentioning another new  proof for $N=2$ in  Mattingly--Sinai
 \cite{Matt99}, which is based on using advanced Fourier Transform
 techniques\footnote{Actually, the proof therein looks
  not that ``elementary" as the title of \cite{Matt99}
 suggests;
  we pretend that our approach is more elementary, though this is not
  a proper point for arguing.}.


\ssk


\subsection{$2m$th-order well-posed Burnett equations}
 \label{SBur}

As a natural extension, similar to KS-type problems in Sections
\ref{SGlH} and \ref{SScal}, we consider the $2m$th-order
 model (\ref{NS1m}).
Actually, for $m \ge 2$, such models are not that formal and have
been known to appear as {\em  Burnett's equations} on the basis of
Grad's method in
 Chapman--Enskog expansions for hydrodynamics.
 Unlike (\ref{NS1m}), the original Burnett equations
 are ill-posed, as backward higher-order parabolic
equations having a wrong sign at the diffusivity operators.
 Namely,   Grad's method applied to kinetic
 equations yields, in addition to the classic operators of the Euler equations, other
  viscosity parts, as follows:
 $$
  \mbox{$
  \vv_t +(\vv \cdot \n)\vv =  
  \sum\limits_{n=0}^\infty \e^{2n+1} \D^n(\mu_n \D \vv)+...= \e\big(\mu_0 \D \vv+ \e^2 \mu_1
  \D^2 \vv+...\big)+...\, ,
   $}
  $$
  where $\e>0$ is essentially the {\em Knudsen number} Kn;  see details in
  Rosenau's  regularization approach,
  \cite{Ros89}. In a full model, truncating such series at $n=0$ leads to the
  Navier--Stokes equations (\ref{NS1}) (with $\mu_0>0$), while $n=1$ is associated with the
    Burnett equations. These are ill-posed since, by expansion, $\mu_1>0$,
    so a backward bi-harmonic flow
    occurs, etc.
 We will refer to (\ref{NS1m}) for $m \ge 2$
as to the {\em well-posed Burnett equations}.

Note also that Burnett-type equations with a small parameter
appeared as higher-order viscosity approximations of the
Navier--Stokes equations is an effective tool for proving
existence of their weak (``turbulent" in Leray's sense) solutions;
see Lions' classic monograph \cite[\S~6, Ch.~1]{JLi}.




 \ssk

  For the system
(\ref{NS1m}), we use the following parameters of scaling:
 \be
 \label{pa1n}
 a_k= C_k^{-\frac pN}, \quad b_k=a_k^{2m}, \andA D_k=C_k^{1+
 \frac{p(2m-1)}N}.
  \ee
This gives the parameter of the convection term in the analogy of
(\ref{NS1n}) as:
 \be
  \label{pa2n}
   \nu_k= C_k^{1- \frac{p(2m-1)}N}.
    \ee
Hence $\nu_k \to 0$ as $k \to \iy$ under the following hypothesis:

\begin{proposition}
 \label{Th.Exm}
 Under the given assumptions,  $L^p$-solutions of the
 well-posed Burnett equations $(\ref{NS1m})$, with
\be
 \label{pa3n}
  \tex{
  p >  p_0= \frac N{2m-1}
  }
  \ee
 in $(\ref{rr1})$,
 do not blow-up, and are uniformly bounded, i.e.,  $\eqref{ka1}$ holds.
 \end{proposition}

The scaling analysis can be applies also directly to the locally
smooth integral nonlocal parabolic flow similar to \eqref{ww1N}.
The critical case $p=p_0$ is then treated by an additional use of
the $(T-t)$-scaling, where some technical difficulties may occur.


\subsection{Well-posed Burnett equations: no blow-up for $N \le
2(2m-1)$}

This  is a simple consequence of the previous scaling analysis. We
recall that, for (\ref{NS1m}), $L^2$-norm of $\vv(t)$ does not
increase with time, so, for smooth solutions, the analogy of
(\ref{ss1}) holds.

\begin{proposition}
 \label{Th.Exm1}
Let $\vv_0 \in L^2(\ren)\cap H^{2m}(\ren)$ be divergence free.
 Then there exists the unique global bounded smooth solution of the
 well-posed Burnett equations $(\ref{NS1m})$ if
\be
 \label{pa3nN}
  \tex{
  N <  2(2m-1).
  }
  \ee
 \end{proposition}

Indeed, substituting $p=2$ into (\ref{pa2n}) yields that

\ssk

\noi(i) $\nu_k \to 0$ in the subcritical case $N < 2(2m-1)$, where
the proof is completed in similar lines and causes no extra
difficulties.

\ssk

 \noi (ii) In addition, $\nu_k \equiv 1$ in the critical case
\be
 \label{w78}
 N=2(2m-1),
  \ee
   where the dynamical system is a gradient one in the
 class of smooth bounded solutions. However, the proof of global
 existence according to the corresponding scaling quite similar to that
 in (\ref{w1}),
 \be
 \label{w1N}
 \tex{
 \vv(x,t)= (T-t)^{-\frac{2m-1}{2m}}\,  \ww(y,\t), \quad y= \frac
 x{(T-t)^{1/2m}}, \quad \t=-\ln(T-t),
  }
  \ee
 will lead to the following rescaled equation:
   \be
   \label{w2N}
    \tex{
     \ww_\t+ {\mathbb P}( \ww \cdot \n)  \ww= \BB^*
    \ww, \quad \BB^*=-(- \D)^m - \frac 1{2m}\, y \cdot \n  -\frac{2m-1}{2m}\,I.
     }
     \ee
Therefore, a proper spectral theory of generalized solenoidal
Hermite polynomials as eigenfunctions of $\BB^*$, with the
spectrum
 \be
 \label{w3N}
  \tex{
  \s(\BB^*)= \big\{ \l_k=-\frac{2m-1}{2m}- \frac k{2m},
  \,\,k=0,1,2,...\big\},
  }
  \ee
  is necessary. This is developed along the lines presented in
  \cite{Eg4}, where the scalar $2m$th-order case was under
  scrutiny.
Eventually, this will allow to treat (we do not guarantee that all
steps are going to be simple as in Section \ref{S6.3}) the global
existence of bounded classical solution in the critical case
(\ref{w78}).
As we have seen, for $m=1$, (\ref{w78}) reads
 $
 N =2
 $
and reflects classic  Leray's
 \cite{Ler33N2}
 and Ladyzhenskaya's
 \cite{Lad58, Lad61, Lad70} results.

As is well-known,  existence or nonexistence of $L^2$-bounded
$L^\iy$-blow-up patterns in the complement to (\ref{pa3nN}),
(\ref{w78}) range
 \be
 \label{ad1}
  N >  2(2m-1), \quad \mbox{or $N  \ge 3$ for $m=1$, including the crucial 3D $N=3$},
   \ee
 comprises the core of the Millennium Prize Problem for the Clay Institute; see Fefferman
 \cite{Feff00}. Concerning possible structures of blow-up
 patterns in the NSEs, see discussion and  review of a large amount of existing literature
  in a \cite{GalNSE} that seems to be a most recent survey in the area.


\end{document}